\documentclass[11pt,a4paper]{article}
\usepackage{amsmath, amssymb}
\usepackage{latexsym}
\usepackage{graphicx, epsfig}
\usepackage{color}
\usepackage{algorithm,algcompatible}

\usepackage{caption}
\usepackage{subcaption}

\newcommand{\DO}{\mathcal{D}_{\partial \Omega}}

\newcommand{\KO}{\mathcal{K}_{\partial \Omega}}

\newcommand{\bef}{\begin{figure}}
\newcommand{\enf}{\end{figure}}
\numberwithin{equation}{section}

\newcommand{\ds}{\displaystyle}

\newcommand{\p}{\partial}

\newcommand{\Rbb}{\mathbb{R}}

\newcommand{\Dcal}{\mathcal{D}}

\newcommand{\Kcal}{\mathcal{K}}

\newcommand{\Scal}{\mathcal{S}}

\newcommand{\Gvf}{\varphi}

\newcommand{\GG}{\Gamma}

\newcommand{\GO}{\Omega}

\newcommand{\beq}{\begin{equation}}
\newcommand{\eeq}{\end{equation}}

\begin{document}

\title{A fast and efficient numerical method for computing the stress concentration between closely located stiff inclusions of general shapes}

\author{Xiaofei Li\thanks{School of mathematical sciences, Zhejiang University of Technology, Hangzhou, 310023, P. R. China (xiaofeilee@hotmail.com, shengqilin@zjut.edu.cn, haojiewang@zjut.edu.cn).}\and
Shengqi Lin\footnotemark[2] \and
Haojie Wang\footnotemark[2]}

\date{\today}
\maketitle

\begin{abstract}

When two stiff inclusions are closely located, the gradient of the solution to the Lam\'{e} system, in other words the stress, may become arbitrarily large as the distance between two inclusions tends to zero. To compute the gradient of the solution in the narrow region, extremely fine meshes are required. It is a challenging problem to numerically compute the stress near the narrow region between two inclusions of general shapes as their distance goes to zero. A recent study \cite{KLY2} has shown that the major singularity of the gradient can be extracted in an explicit way for two general shaped inclusions. Thus the complexity of the computation can be greatly reduced by removing the singular term and it suffices to compute the residual term only using regular meshes. The goal of this paper is to numerically compute the stress concentration in a fast and efficient way. In this paper, we compute the value of the stress concentration factor, which is the normalized magnitude of the stress concentration, for general shaped domain as the distance between two inclusions tends to zero. We also compute the solution for two closely located inclusions of general shapes and show the convergence of the solution. Only regular meshes are used in our numerical computation and the results clearly show that the characterization of the singular term method can be efficiently used for computation of the stress concentration between two closely located inclusions of general shapes.

\end{abstract}


%%%%%%%%%%%%%%%%%%%%%%%%%%%%%%%%%%%%%%%%%%%%%%%%%%%%%%%%%%%%
\section{Introduction}
%%%%%%%%%%%%%%%%%%%%%%%%%%%%%%%%%%%%%%%%%%%%%%%%%%%%%%%%%%%%

\begin{figure}[!ht]
\centering
\epsfig{figure=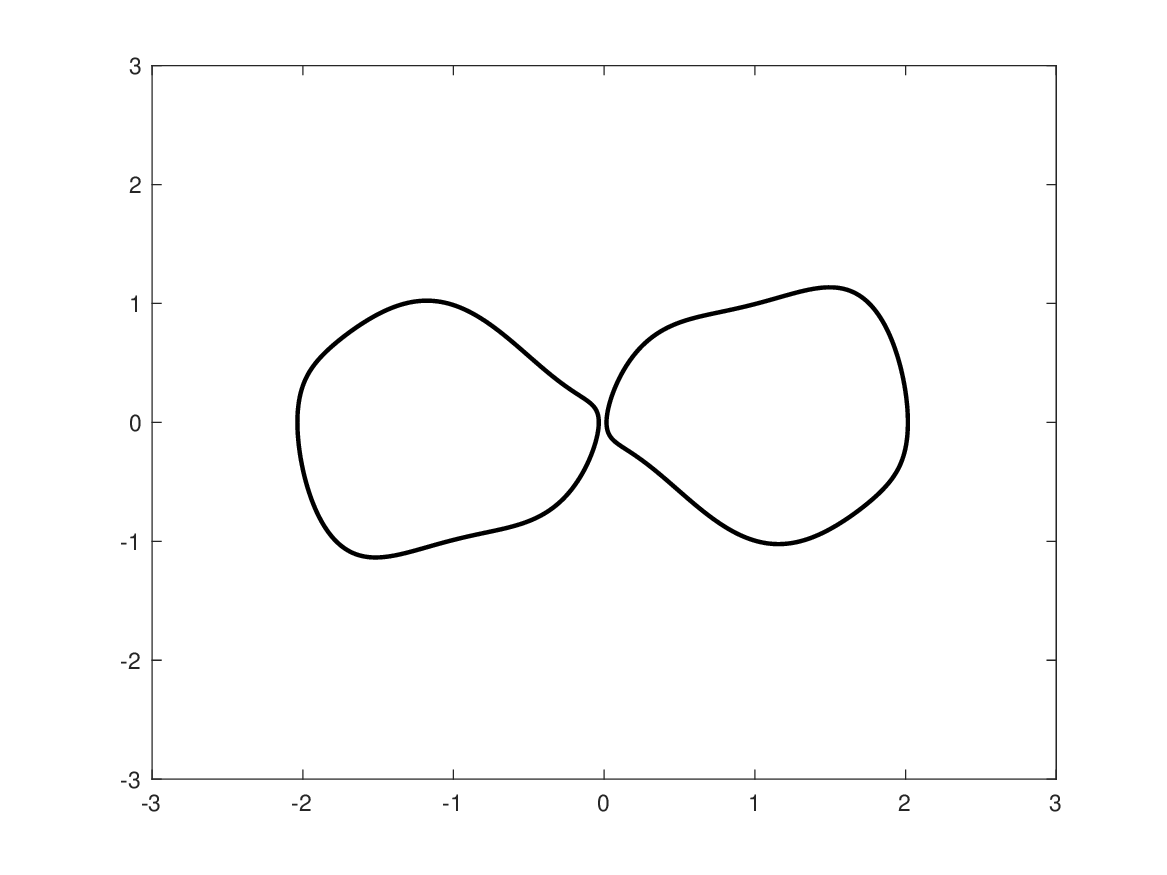, height=5cm}
\caption{General geometry.}
\label{geom}
\end{figure}

Let $D_1$ and $D_2$ be two closely located strictly convex simply connected domains in $\mathbb{R}^2$ with $C^{2,\gamma}$ smooth boundaries for some $\gamma \in (0,1)$, see, for example, Figure \ref{geom}. Let
$$\epsilon: = \mbox{dist}(D_1,D_2),$$
which is assumed to be small. 
Assume that there are unique points $z_1\in \p D_1$ and $z_2\in \p D_2$ such that
$$|z_1 - z_2| = \mbox{dist}(D_1,D_2),$$
where $z_1\in \p D_1$ and $z_2\in \p D_2$ are the closest points. One can further relax the (global) strict convexity assumption of $D_1$ and $D_2$ by assuming that $D_j$ is strictly convex near $z_j$, $j = 1,2$, namely, there is a common neighborhood $U$ of $z_1$ and $z_2$ such that $D_j \cap U$ is strictly convex for $j = 1,2$. Moreover, we assume that
$$\mbox{dist}(D_1, D_2\backslash U) \geq C \quad \mbox{and} \quad \mbox{dist}(D_2, D_1\backslash U) \geq C $$
for some positive constants $C$ independent of $\epsilon$. Note that strictly convex domains satisfy all the assumptions.

Let $H$ be a given entire harmonic function in $\mathbb{R}^2$. We consider the following two dimensional problem
\beq\label{main}
\begin{cases}
\Delta u = 0 & \quad \mbox{in } \Rbb^2 \backslash \overline{(D_{1}\cup D_2)},\\
u  = \lambda_j (\mbox{constant}) &\quad \mbox{on } \partial D_{j}, ~j = 1,2,\\
u -H(x)= O(|x|^{-1}) &\quad \mbox{as } |x| \to \infty,
\end{cases}
\eeq
where the constants $\lambda_j$ are determined by the conditions
$$
\int_{\p D_j} \frac{\p u_\epsilon}{\p\nu}\Big|_{+}ds = 0,~j = 1,2.
$$
Here and throughout this paper, $\nu$ denotes the unit outward normal on $\p D_j$. The notations $|_+$ and $|_-$ are for limits from outside and inside of inclusions, respectively. It is worth mentioning that the constants $\lambda_1$ and $\lambda_2$ may or may not be the same depending on the applied field $H$. 
Two inclusions $D_1$ and $D_2$ may represent the two dimensional cross-sections of two parallel elastic fibers embedded in an infinite elastic matrix. The solution $u$ represents the out-of-plane elastic displacement, and the gradient of the solution is proportional to the shear stress. When the inclusions are fiber-reinforced composites that are densely packed, the stress concentration may occur and cause material failure due to the damage of fiber composites. Thus, it is important to quantitatively understand the stress concentration. This problem was first raised in \cite{Babuska}. During the last two decades or more, significant development on this problem has been developed. It has been proved that the gradient blow-up rate is $\epsilon^{-1/2}$ in two dimensions \cite{AKLLL, AKLLZ, AKL, BLL1, BLY, BC, BV, Keller, LY, Yun, Yun2}, and $|\epsilon \ln \epsilon|^{-1}$ in three dimensions \cite{BLL2, BLY, BLY2, KLY3, Lekner, Lekner2, Lekner3, LY2}, see \cite{K} for more references.

As for the numerical computation, due to the stress concentration near the narrow region between two inclusions, extremely fine meshes are required to numerically compute the stress. Recently, a better understanding of the stress concentration has been proposed in \cite{Yun} and used for various different circumstances \cite{ACKLY,KLY1,KLY2,KLY3, KY}. It is shown that the asymptotic behavior of the gradient of the solution can be characterized by a singular function associated with two inclusions as the distance $\epsilon$ tends to zero. Using this singular function, the solution can be decomposed into a singular and a regular term. After extracting the singular term in an explicit way, it is sufficient to compute the residual term only using regular meshes. Therefore, the complexity of computation can be greatly reduced by removing the singular term. In fact, this idea was exploited numerically in \cite{KLY1} for the special case when two inclusions are of circular shapes. The stress concentration factor, which is the normalized magnitude of the stress concentration, is only related to the radii of two disks as well as the applied field which is easy to compute. Meanwhile, the singular function which characterizes the singular behavior of the stress is also explicit for circular domains. The numerical results in \cite{KLY1} showed that the characterization of the singular term method can be effectively used for computation of the gradient in the presence of two nearly touching disks.

However, if two nearly touching inclusions are of general shapes, it is difficult to find accurate stress concentration factor since it is related to the global information of the singular function which is not explicit for general shapes. Thanks to the study in \cite{KLY2}, where it shows that as the distance between two inclusions tends to zero, the stress concentration factor converges to a certain integral of the solution to the system when two inclusions are touching. Our goal of this paper is to numerically compute the stress concentration between two closely located stiff inclusions of general shapes in a fast and accurate way. Motivated by the idea of characterization of the singular term and the fact that the stress concentration factor converges to a certain integral of the solution to the touching case which is shown in \cite{KLY2}, we are able to obtain fast and accurate numerical computation of the stress concentration only using regular meshes. We will compute the value of the stress concentration factor accurately by numerically solving the touching case. We then numerically compute the solution for two closely located inclusions of general shapes and show the convergence of the solution. Our numerical computations clearly show that the characterization of the singular term method can be efficiently used for computation of the stress concentration when two closely located inclusions are of general shapes. It is worth mentioning that only regular meshes are used in our computations.

This paper is organized as follows. In section 2, we briefly review on the characterization of the singular term method. In section 3, we show how to compute the stress concentration factor by solving the touching problem. In section 4, we give the numerical computation scheme of the characterization of singular term method and show the effectiveness of it. Some numerical examples of the computation of the solution for  general shaped domains are given in Section 5. This paper ends with a short inclusion.

%%%%%%%%%%%%%%%%%%%%%%%%%%%%%%%%%%%%%%%%%%%%%%%%%%%%%%%%%%%%
\section{Review on the characterization of the singular term method}
%%%%%%%%%%%%%%%%%%%%%%%%%%%%%%%%%%%%%%%%%%%%%%%%%%%%%%%%%%%%

In this section, we briefly review the characterization of the singular term method in \cite{ACKLY, KLY2}. 

Let $D_1$ and $D_2$ be two stiff inclusions in $\mathbb{R}^2$ with $C^{2,\gamma}$ smooth boundaries for some $\gamma \in (0,1)$. They satisfy the geometric description in the previous section. Let $z_1\in \p D_1$ and $z_2\in \p D_2$ are the closest points such that $|z_1 - z_2| = \mbox{dist}(D_1,D_2)$. Let $\epsilon: = \mbox{dist}(D_1,D_2)$. After rotation and translation, we assume that $z=(z_1+z_2)/2$ is at the origin. We also assume that the $x_1$-axis is parallel to the vector $z_2-z_1$. Then
$$z_1 = (-\epsilon/2,0) \quad \mbox{and} \quad z_2 = (\epsilon/2,0).$$

Let $B_j$ be the disk osculating to $D_j$ at  $z_j$ $(j = 1,2)$. Let $R_j$ be the reflection with respect to $\p B_j$, $j = 1,2$, and let $p_1 \in B_1$ is the fix point of the mixed reflection $R_1R_2$, and $p_2\in D_2$ be that of $R_2R_1$.  Let $q$ be the singular function associated with $B_1$ and $B_2$, given as follows:
\beq\label{qB}
q (x) = \frac{1}{2\pi} \left(\ln |x-p_1| - \ln|x-p_2| \right). 
\eeq
It is easy to see that $\nabla q$ blows up at the order of $\epsilon^{-1/2}$ near the narrow region between two inclusions. 

It is proved in \cite{ACKLY, KLY2} that the solution $u$ to the problem \eqref{main} admits the following representation:
\beq\label{blow}
\nabla u(x) = \alpha_0 \nabla q(x) (1+O(\epsilon^{\gamma/2})) + O(1) \quad \mbox{as} ~\epsilon\rightarrow 0.
\eeq
Since $\nabla q$ blows up at the order of $\epsilon^{-1/2}$, $\alpha_0 \nabla q$ is the singular part of $\nabla u$. Here, $\alpha_0$ is the so called stress concentration factor. In particularly, if two inclusions are disks of radii $r_1$ and $r_2$, respectively, the stress concentration factor is given as
$$\alpha_0 =  \frac{4\pi r_1 r_2}{r_1+r_2} (\nu \cdot \nabla H)(z).$$

For general shaped inclusions, the stress concentration factor was shown, in  \cite{KLY2}, to converge to a certain integral of the solution to the touching case as the distance between two inclusions tends to zero, namely, the case when $\epsilon =0$. In fact, it is shown in \cite{KLY2} that $\alpha_0$ can be computed in the following way. For $\rho>0$, let
\beq \label{Drho}
D_{\rho} = (D_1^0 \cup D_2^0) \cup ([-\rho,\rho]\times [-\rho,\rho]),
\eeq
which is of dumbbell shape. Let $u_\rho$ be the solution to
\beq\label{u-rho}
\begin{cases}
\Delta u_\rho = 0 \quad \mbox{in } \Rbb^d \backslash \overline{D}_{\rho},\\
u_\rho = \lambda_\rho (\mbox{constant}) \quad \mbox{on } \partial D_{\rho},\\
u_\rho(x)-H( x)= O(|x|^{-1}) \quad \quad\mbox{as } |x| \to \infty,
\end{cases}
\eeq 
where the constant $\lambda_\rho$ is determined by the additional condition
$$
\int_{\p D_\rho} \p_\nu u_\rho|_+  ds=0.
$$
Let
\beq\label{alpha-rho}
\alpha_\rho = \int_{\p D_1^0 \backslash [-2\rho,2\rho]\times [-2\rho,2\rho]} \p_\nu u_\rho|_+ ds.
\eeq
Then there are constants $C$ and $A>0$ independent of $\rho$ such that
\beq\label{equiv}
|\alpha_0 - \alpha_\rho|\leq C \exp \left(-\frac{A}{\rho}\right).
\eeq
By \eqref{equiv}, one can obtain an accurate approximation of the stress concentration factor $\alpha_0$ by computing \eqref{alpha-rho} through the touching case \eqref{u-rho}. Therefore, by \eqref{blow}, the solution $u$ to the problem \eqref{main} can be decomposed as a singular term and a regular term:
$$u(x) = \alpha_0 q(x) + b(x),$$
where $q$ is given by \eqref{qB} and 
$$\|\nabla b\|_{L^{\infty}  (\mathbb{R}^2\backslash D_1\cup D_2)} \leq C$$
for a constant $C$ independent of $\epsilon$. For the numerical computation, it is sufficient to compute $b$ only using regular meshes. As it is shown above, it is a key point to obtain accurate value of the stress concentration factor $\alpha_0$. In the next section, we will show how to obtain the accurate value of the stress concentration factor through numerically solving the touching case.

%%%%%%%%%%%%%%%%%%%%%%%%%%%%%%%%%%%%%%%%%%%%%%%%%%%%%%%%%%%%
\section{Computation of the stress concentration factor $\alpha_0$}
%%%%%%%%%%%%%%%%%%%%%%%%%%%%%%%%%%%%%%%%%%%%%%%%%%%%%%%%%%%%

In this section, we compute the stress concentration factor $\alpha_0$ by solving the touching case \eqref{u-rho} using boundary element method by Matlab. We also show the convergent rate of the computation. Before doing so, we introduce some basic concepts on layer potentials.  

Let
\begin{equation*}
\GG(x)=\frac{1}{2\pi} \ln |x|,
\end{equation*}
the fundamental solution to the Laplacian in two dimensions. Let $\GO$ be a simply connected domain with the Lipschitz boundary. The single and double layer potentials of a function $\Gvf$ on $\p \GO$ are defined to be
\begin{align*}
\Scal_{\p \GO}[\Gvf](x) &:= \int_{\p \GO} \GG(x-y) \Gvf(y)\, ds(y), \quad x\in \Rbb^2, \\
\DO[\Gvf](x) &:= \int_{\p \GO} \p_{\nu_y}\GG(x-y) \Gvf(y)\, ds(y), \quad x\in \Rbb^2 \setminus \p \GO,
\end{align*}
where $\p_{\nu_y}$ denotes outward normal derivative with respect to $y$-variables. It is well known (see, for example, \cite{AmKa07Book2}) that the single and double layer potentials satisfy the following jump relations:
\begin{align}
\p_\nu \Scal_{\p \GO}[\Gvf](x)\Big|_{\pm} &=(\pm \ds\frac{1}{2}I+\Kcal_{\p \GO}^{*})[\Gvf](x), \quad\mbox{a.e. } x\in\p \GO, \label{singlejump} \\
\Dcal_{\p \GO}[\Gvf](x)\big|_{\pm} &=(\mp\ds\frac{1}{2}I+\Kcal_{\p \GO})[\Gvf](x), \quad\mbox{a.e. } x\in\p \GO,
\label{doublejump}
\end{align}
where the operator $\KO$ on $\p \GO$ is defined by
\begin{equation*}
\Kcal_{\p \GO}[\Gvf](x)= \mbox{p.v.} \int_{\p \GO} \p_{\nu_y}\GG(x-y)  \Gvf(y)\, ds(y),
\end{equation*}
and $\Kcal_{\p \GO}^*$ is the $L^2$-adjoint of $\Kcal_{\p\GO}$. Here, p.v. stands for the Cauchy principal value.

The stress concentration factor can be precisely estimated by the integral \eqref{alpha-rho} through solving the touching problem \eqref{u-rho}, which is demonstrated in the previous section. By layer potential techniques, the solution $u_\rho$ to \eqref{u-rho} can be represented as
\beq\label{u_rho}
u_\rho(x) = H(x) + \mathcal{S}_{\p D_\rho}[\psi](x),\quad x\in \mathbb{R}^2 \backslash \overline{D}_\rho,
\eeq
for $\psi \in L^2_0(\p D_\rho)$, where $L^2_0$ denotes the set of $L^2$ functions with mean zero. Since $u_\rho$ is constant on $\p D_\rho$, $\psi$ should satisfy
$$\frac{\p H}{\p\nu} + \frac{\p}{\p \nu} \Scal_{\p D_\rho}[\psi] \Big|_{-} = 0   \quad \mbox{on}~\p D_\rho,$$
which, according to \eqref{singlejump}, can be written as 
\beq\label{psi}
\frac{\p H}{\p\nu} + \left(-\frac{1}{2} I + \mathcal{K}^*_{\p D_\rho}  \right) [\psi] = 0  \quad \mbox{on}~\p D_\rho.
\eeq
Taking outward normal derivative of \eqref{u_rho}, and by the jump relation \eqref{singlejump}, we have
\beq\label{urho-de}
\frac{\p u_\rho}{\p\nu} \Big|_+ = \frac{\p H}{\p\nu} + \left(\frac{1}{2} I + \mathcal{K}^*_{\p D_\rho}  \right) [\psi]  \quad \mbox{on}~\p D_\rho.
\eeq
In view of \eqref{psi}, then \eqref{urho-de} becomes
\beq\label{urho-nu}
\frac{\p u_\rho}{\p\nu} \Big|_+ = \psi  \quad \mbox{on}~\p D_\rho.
\eeq
Hence by \eqref{alpha-rho} and \eqref{urho-nu}, we have
\beq\label{factor}
\alpha_\rho =  \int_{\p D_1^0 \backslash [-2\rho,2\rho]\times [-2\rho,2\rho]} \psi ds,
\eeq
where $\psi$ is given by \eqref{psi} as follows
$$\psi = \left(\frac{1}{2} I - \mathcal{K}^*_{\p D_\rho}  \right)^{-1}\left[\frac{\p H}{\p\nu}\right ].$$
The density function $\psi$ can be uniquely solved. In fact, it is well known, see for example \cite{AmKa07Book2}, that the operator $\lambda I - \mathcal{K}^*_{\p D_\rho}$ is one to one on $L^2_0(\p D_\rho)$ if $|\lambda| \geq 1/2$ when $D_\rho$ is a bounded Lipschitz domain. 

\begin{figure}[!ht]
\centering
\epsfig{figure=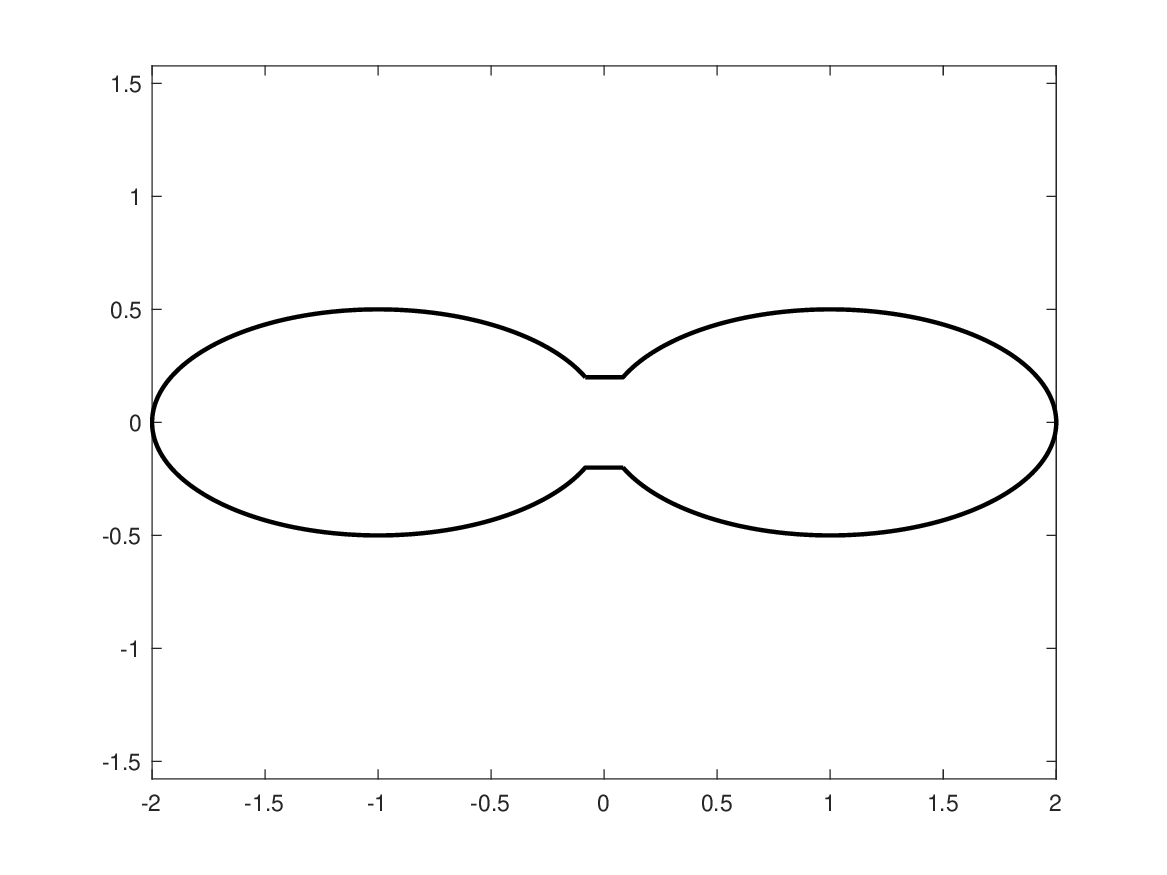, height=4cm}
\caption{$D_\rho$ for two touching ellipses.}
\label{neck}
\end{figure}

We now compute \eqref{factor} using boundary element method. For example, let $D_1$ and $D_2$ be two elliptic inclusions of the same major axis $a = 2$ and minor axis $b = 1$, centered at $(-a-\epsilon/2,0)$ and $(a+\epsilon/2,0)$, respectively, where $\epsilon = 0.01$. The domain $D_\rho$ which is defined by \eqref{Drho} is the dumbbell shaped domain shown in Figure \ref{neck}.
Discretize each boundary of $D_j$, $j=1,2$ into $N$ points and each connecting segment between $D_1$  and $D_2$ into $N/16$ points. 

\begin{figure}[!ht] 
\centering
\epsfig{figure=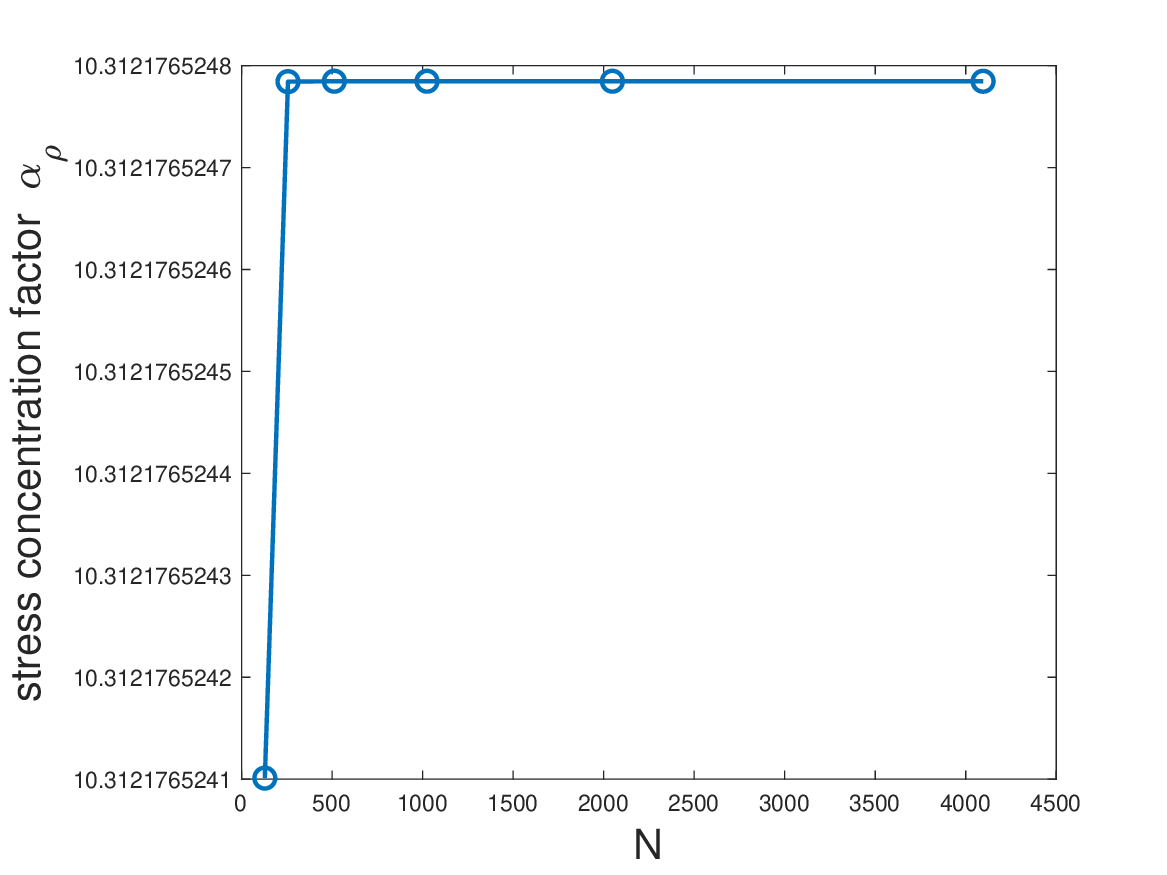, height=3.5cm}\epsfig{figure=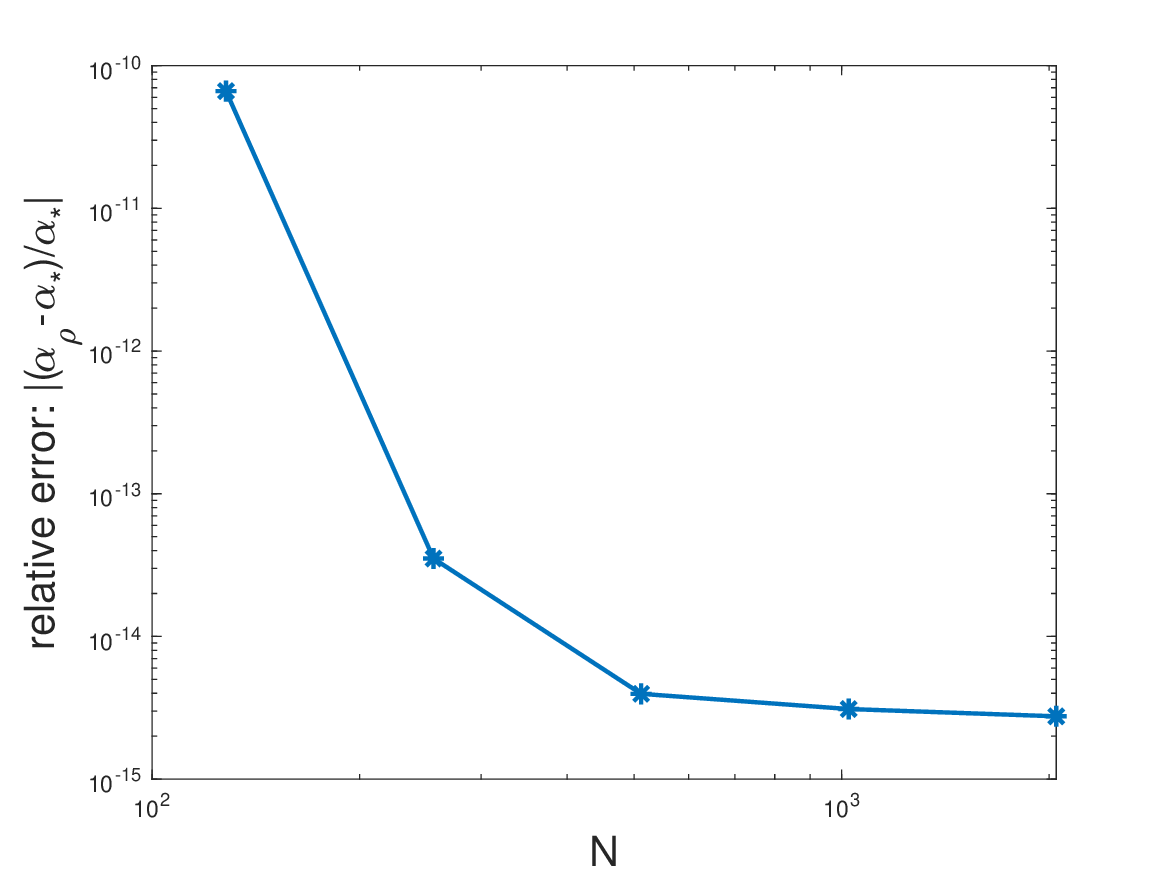, height=3.5cm}\epsfig{figure=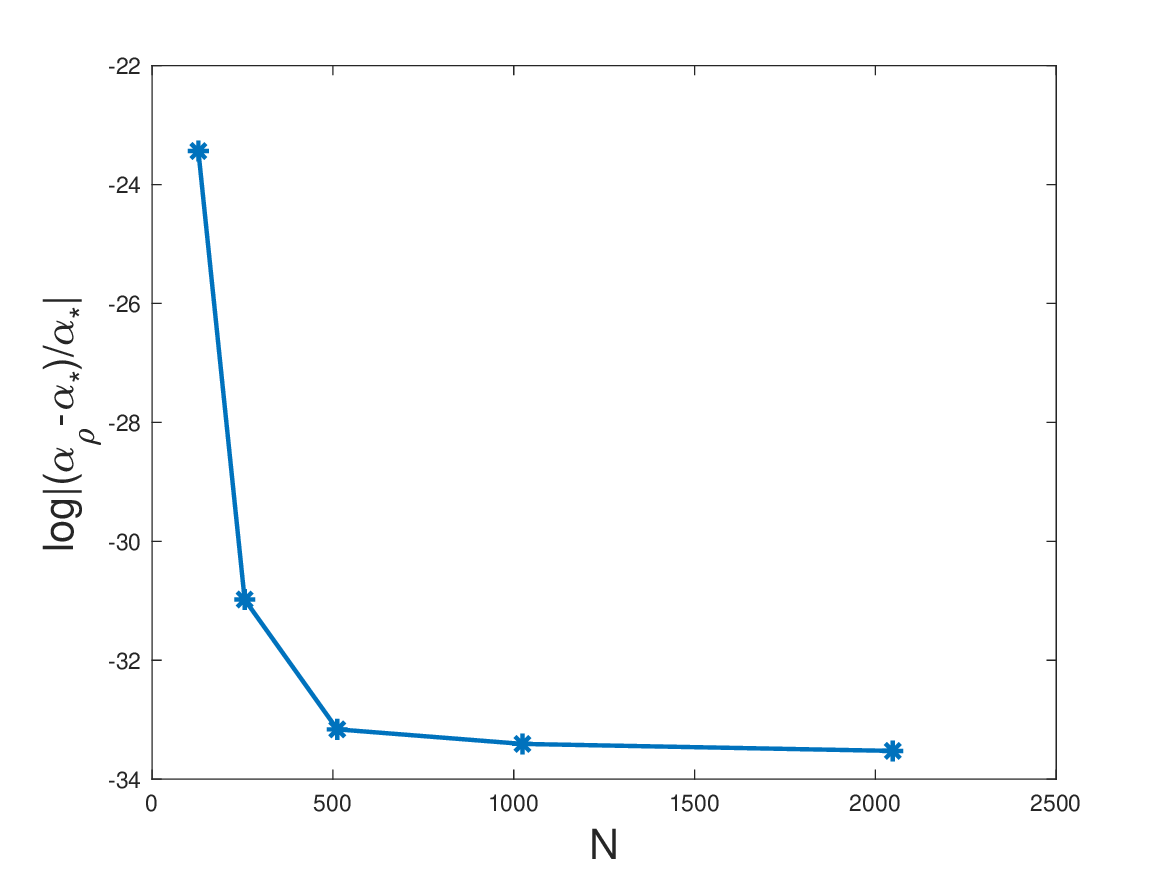, height=3.5cm}
\caption{Left: $\alpha_\rho$ for different values of the grids number $N$. Middle: The relative error: $|(\alpha_\rho-\alpha_*)/\alpha_*|$. Right: The convergent rate: $\log|(\alpha_\rho-\alpha_*)/\alpha_*|$.}
\label{alpha_N}
\end{figure}

Firstly, we fix $\rho = 0.05$ in \eqref{factor} and change the values of the grids number $N$ with $N = 256,512,1024,2048,4096$. Figure \ref{alpha_N} shows the numerical result of $\alpha_\rho$ for different values of $N$. One can easily see that $\alpha_\rho$ converges as the number of grids points increases. Denote $\alpha_*$ as the value of $\alpha_{\rho} $ with finer grids $N= 4096$. We then compare each $\alpha_\rho$ with $\alpha_*$. The relative error is shown in Figure \ref{alpha_N} (Middle). The convergent rate is shown in Figure \ref{alpha_N} (Right). One can see that $\alpha_\rho$ converges very fast as $N$ increases.

\begin{figure}[!ht] 
\centering
\epsfig{figure=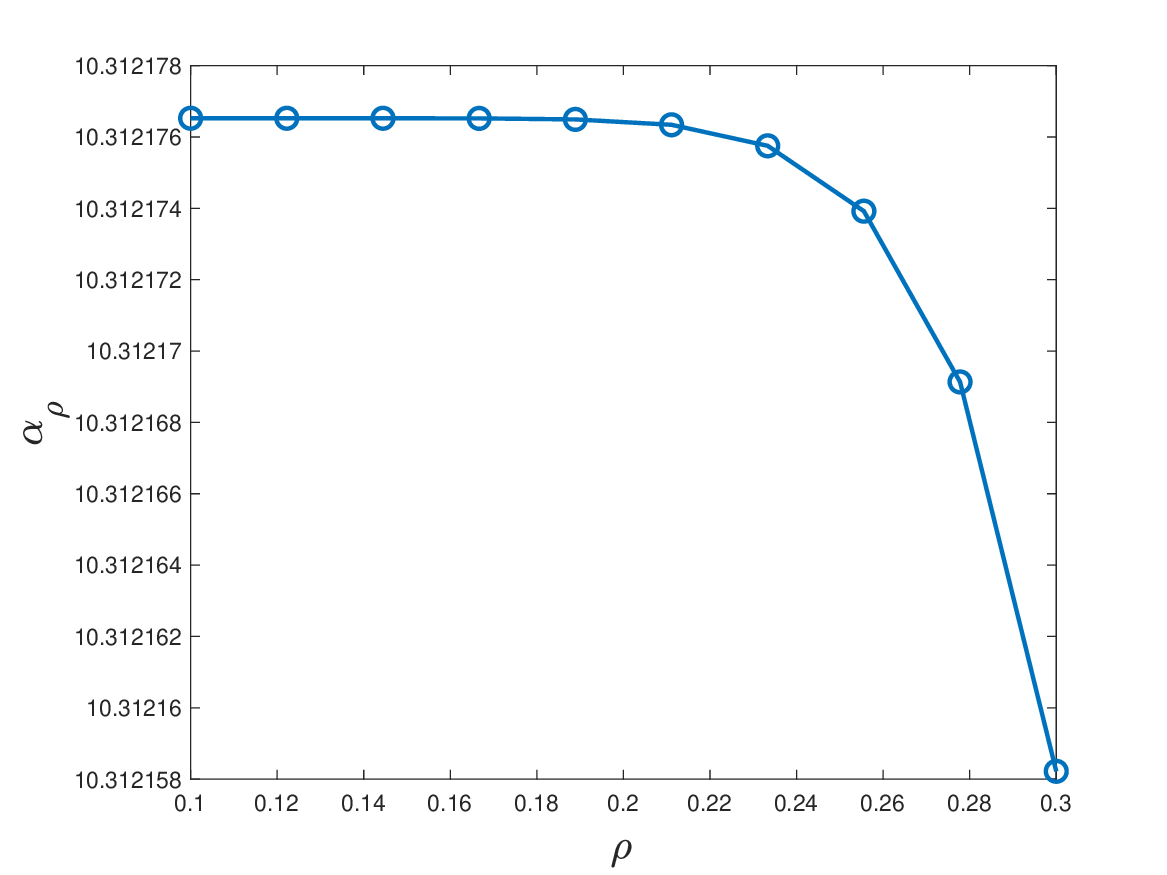, height=3.5cm}\epsfig{figure=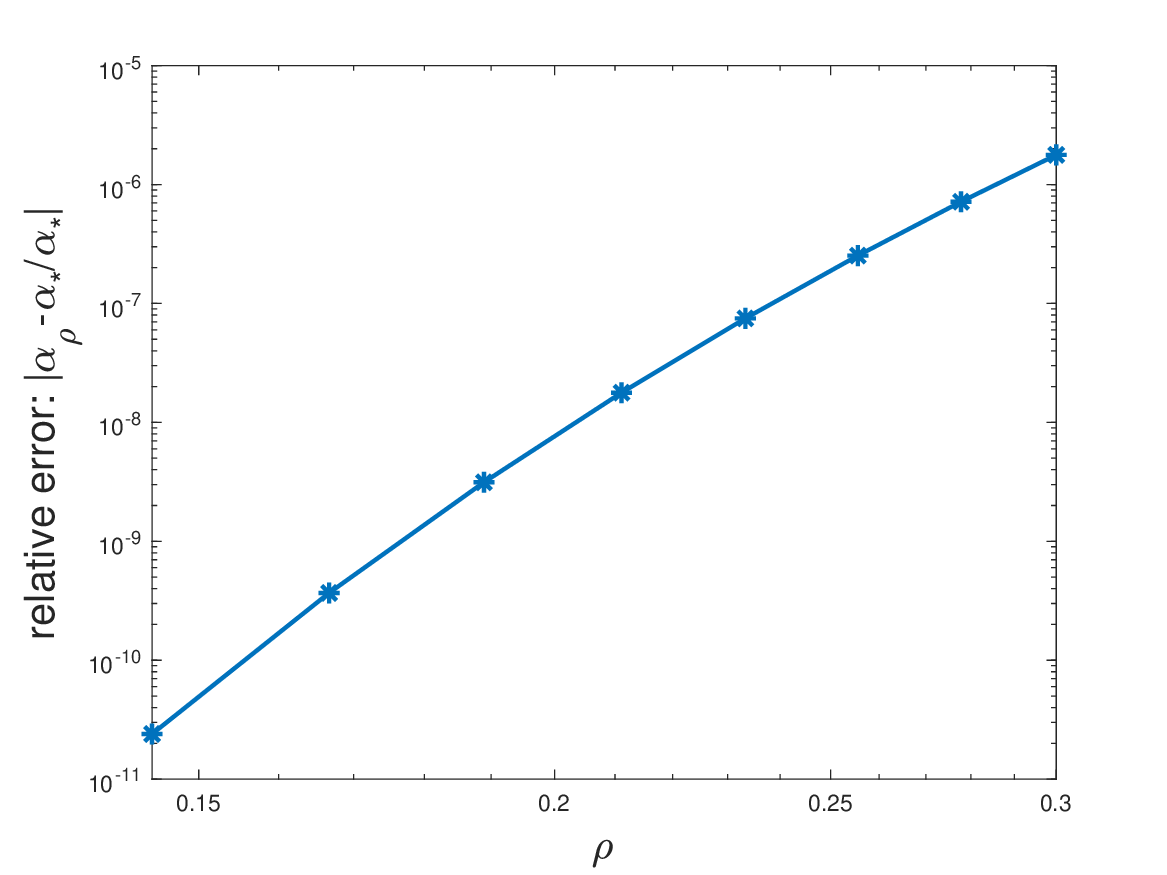, height=3.5cm}\epsfig{figure=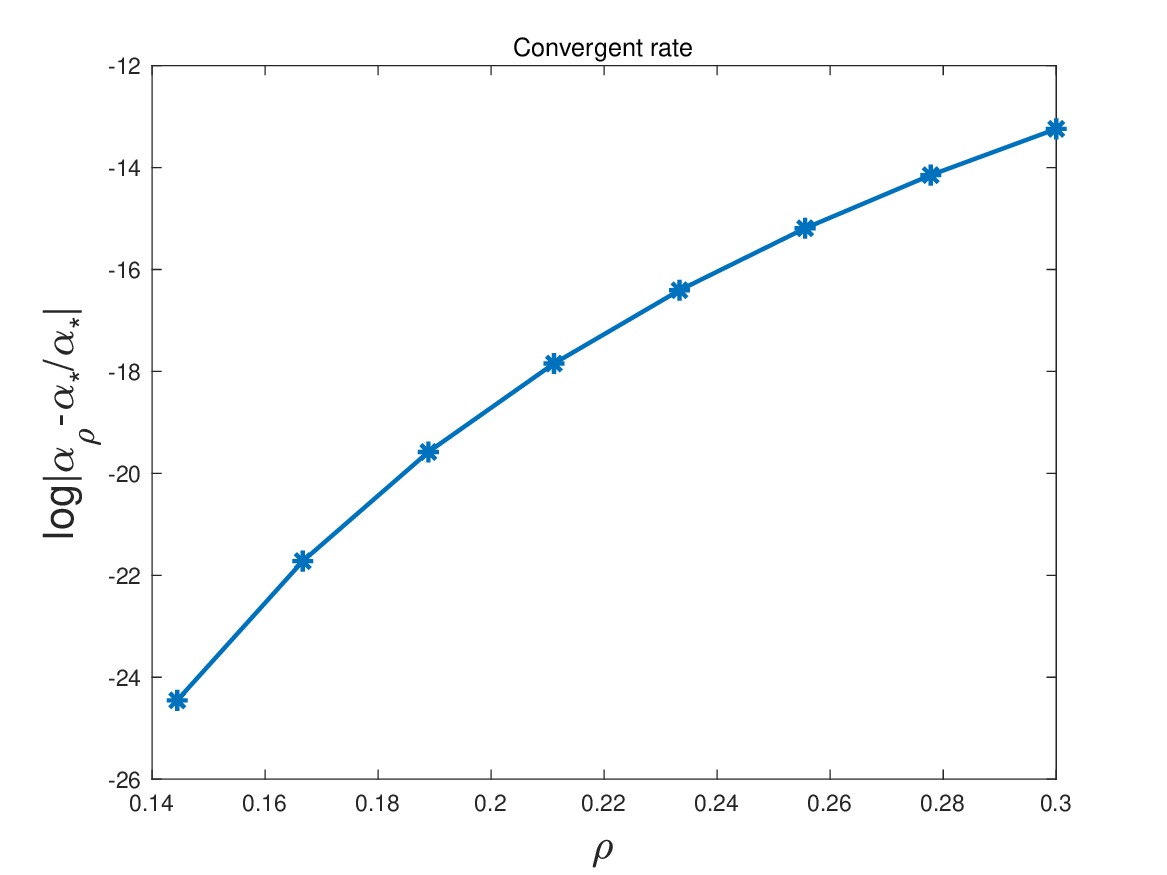, height=3.5cm}
\caption{Left: $\alpha_\rho$ for different values of $\rho$. Middle: The relative error: $|(\alpha_\rho-\alpha_*)/\alpha_*|$. Right: The convergent rate: $\log|(\alpha_\rho-\alpha_*)/\alpha_*|$.}
\label{alpha_rho}
\end{figure}

Secondly, we fix $N = 512$ and change $\rho$ from $0.3$ to $0.1$. Figure \ref{alpha_rho} shows the numerical values of $\alpha_\rho$ for different values of $\rho$. From the left-hand side figure one can see that $\alpha_\rho$ converges as $\rho$ decreases. Denote $\alpha_*$ as the value of $\alpha_\rho$ when $\rho= 0.1$. Then the relative errors of $\alpha_\rho$ and the convergent rate is shown in the middle and right-hand side figure of Figure \ref{alpha_rho}, respectively. One can also see that $\alpha_\rho$ converges very fast as $\rho$ decreases.

Figure \ref{alpha_N} and Figure \ref{alpha_rho} both show that we can obtain accurate value of the stress concentration factor by numerical computation.

%%%%%%%%%%%%%%%%%%%%%%%%%%%%%%%%%%%%%%%%%%%%%
\section{Numerical computations}
%%%%%%%%%%%%%%%%%%%%%%%%%%%%%%%%%%%%%%%%%%%%%

In this section, we provide numerical scheme on the computation of the solution to \eqref{main} using characterization of the singular term method. We show that it can be efficiently used for the computation of the stress concentration by comparing the convergent rate with the solution computed using layer potential techniques in a direct way. The  boundary element method is used for both methods. 

Before providing the numerical scheme, we first derive the related system of integral equations.

\subsection{System of the integral equations}
%%%%%%%%%%%%%%%%%%%%%%%%%%%%%%%%%%%%%%%%%%%%%

Let $D_1$ and $D_2$ be the same as in the previous section. Let $B_j$ be the disk osculating to $D_j$ at $z_j$ ($j = 1,2$), where $z_1 = (-\epsilon/2,0)$ and $z_2 = (\epsilon/2,0)$. Let $\kappa_j$ be the curvature of $D_j$ at $z_j$. Then the radius of $B_j$ is $r_j = 1/\kappa_j$, $j = 1,2$. Let $c_j$ be the center point of the disk $B_j$, $j = 1,2$.

Define singular function $q$ in the spirit of \eqref{qB} as follows
\beq\label{q}
q (x) = \frac{1}{2\pi} \left(\ln |x-p_1| - \ln|x-p_2| - \ln|x-c_1| + \ln|x-c_2| \right),
\eeq
for $x\in \mathbb{R}^2\backslash \overline{(D_{1}\cup D_2)}$, where $p_1$ and $p_2$ are two fix points of the mixed reflection with respect to $\p B_j$, $j = 1,2$. In fact, it is shown in \cite{Yun, Yun2} that the fixed points $p_1$ and $p_2$ are given by 
\beq\label{p12}
p_1 = \left(-\sqrt{2} \sqrt{\frac{r_1 r_2}{r_1+r_2}} \sqrt{\epsilon} + O(\epsilon), 0 \right) \quad \mbox{and} \quad p_2 = \left(\sqrt{2} \sqrt{\frac{r_1 r_2}{r_1+r_2}} \sqrt{\epsilon} + O(\epsilon), 0 \right).
\eeq

In view of \eqref{blow}, we look for a solution $u$ to \eqref{main} in the following form
\beq\label{u-eps}
u(x) = \alpha_0 q (x) + H(x) + \mathcal{S}_{\p D_1}[\phi_1](x) + \mathcal{S}_{\p D_2}[\phi_2](x), \quad x\in \Rbb^2 \backslash \overline{(D_{1}\cup D_2)},
\eeq
where $(\phi_1,\phi_2) \in L^2_0(\p D_1)\times L^2_0(\p D_2)$ are to be determined. It is worth mentioning that the gradient of $H+\mathcal{S}_{\p D_1}[\phi_1] + \mathcal{S}_{\p D_2}[\phi_2]$ is bounded on $\Rbb^2 \backslash \overline{(D_{1}\cup D_2)}$  according to \eqref{blow}, and hence $\|\phi_1\|_{L^\infty(\p D_1)}$ and $\|\phi_2\|_{L^\infty(\p D_2)}$ are bounded regardless of $\epsilon$. We use the fact that $\frac{\partial u}{\partial \nu}|_{-} = 0$ on $\partial D_j$, $j = 1,2$ to find the integral equations for $(\phi_1,\phi_2)$. In order to do so, we take harmonic extension of $u$ toward the interior of $D_1\cup D_2$. Note that $H$, $\mathcal{S}_{\p D_1}[\phi_1]$ and $\mathcal{S}_{\p D_2}[\phi_2]$ are continuous in $\mathbb{R}^2$ and harmonic in $D_1\cup D_2$. Hence, it remains to find the harmonic extension of $q$ toward the interior of $D_1\cup D_2$.

Let $q_j$ be the harmonic extension of $q$ towards the interior of $D_j$, $j = 1,2$, respectively. Then $q_j$ should satisfy the following Dirichlet problem:
\beq\label{qj}
\begin{cases}
\Delta q_j (x) = 0 \quad \mbox{in}~ D_j,\\
q_j(x) = q_j |_{\p D_j},
\end{cases}
\eeq
where the boundary data is given explicitly by \eqref{q}:
$$ \displaystyle  q_j |_{\p D_j} = \frac{1}{2\pi} \left(\ln |x-p_1| - \ln|x-p_2| - \ln|x-c_1| + \ln|x-c_2| \right) \quad \mbox{on} ~\p D_j.$$
By numerically solving \eqref{qj} for each $j = 1,2$, one can obtain the interior harmonic extension of the singular function $q_j$ in $D_j$. Let
\begin{equation*}
q^G (x) =\begin{cases}
q_1 (x) &~ \mbox{in}~ D_1,\\
q_2 (x)&~ \mbox{in}~ D_2,\\
q (x)&~ \mbox{in}~ \Rbb^2 \backslash \overline{(D_{1}\cup D_2)}.
\end{cases}
\end{equation*}
Then $q^G$ is continuous in $\Rbb^2$ and harmonic in $D_1$ and $D_2$ as well as in $\Rbb^2 \backslash \overline{(D_{1}\cup D_2)}$.

Define
$$u^G(x) = \alpha_0 q^G (x) + H(x) + \mathcal{S}_{\p D_1}[\phi_1](x) + \mathcal{S}_{\p D_2}[\phi_2](x), \quad x\in \Rbb^2.$$
Then $u^G$ is continuous in $\mathbb{R}^2$ and harmonic in $D_1$, $D_2$ and $\Rbb^2 \backslash \overline{(D_{1}\cup D_2)}$. Since $u^G$ is constant on $\p D_j$, $j = 1,2$,  $u^G$ should be constant in $D_j$, $j = 1,2$. Taking inward normal derivative of $u^G$ on $\p D_j$ and by the jump relation \eqref{singlejump}, we obtain the following integral equations for $(\phi_1,\phi_2)$:
\beq \label{sys-phi}
\begin{cases}
\displaystyle \left(-\frac{1}{2}I + \mathcal{K}_{\p D_1}^* \right) [\phi_1] + \frac{\p }{\p \nu_{D_1}} \mathcal{S}_{\p D_2}[\phi_2] = -\frac{\p H}{\p \nu_{D_1}} - \alpha_0 \frac{\p q_1}{\p \nu_{D_1}}\Big|_{-} &~ \mbox{on}~\p D_1,\\[1em]
\displaystyle \frac{\p }{\p \nu_{D_2}} \mathcal{S}_{\p D_1}[\phi_1] + \left(-\frac{1}{2}I + \mathcal{K}_{\p D_2}^* \right) [\phi_2] = -\frac{\p H}{\p \nu_{D_2}} - \alpha_0 \frac{\p q_2}{\p \nu_{D_2}}\Big|_{-}  &~ \mbox{on}~\p D_2, 
\end{cases}
\eeq
where $\frac{\p q_j}{\p \nu_{D_j}}\big|_{-}$, $j = 1,2$, can be obtained by solving \eqref{qj} numerically. The density functions $(\phi_1,\phi_2) \in L^2_0(\p D_1)\times L^2_0(\p D_2)$ can be uniquely determined by solving the system of integral equations \eqref{sys-phi}. In fact, denote 
$$\mathbb{I} = \begin{bmatrix} I & 0 \\ 0 & I \end{bmatrix}, \quad  \mathbb{K}^* = \begin{bmatrix} \mathcal{K}_{\p D_1}^* & \frac{\p }{\p \nu_{D_1}} \mathcal{S}_{\p D_2} \\ \frac{\p }{\p \nu_{D_2}} \mathcal{S}_{\p D_1} & \mathcal{K}_{\p D_2}^* \end{bmatrix},$$
then $-\frac{1}{2} \mathbb{I} + \mathbb{K}^*$ is invertible on $L^2_0(\p D_1)\times L^2_0(\p D_2)$, which is shown in, for example, \cite{ACKLY}.

\begin{figure}[!ht]
\begin{center}
\epsfig{figure=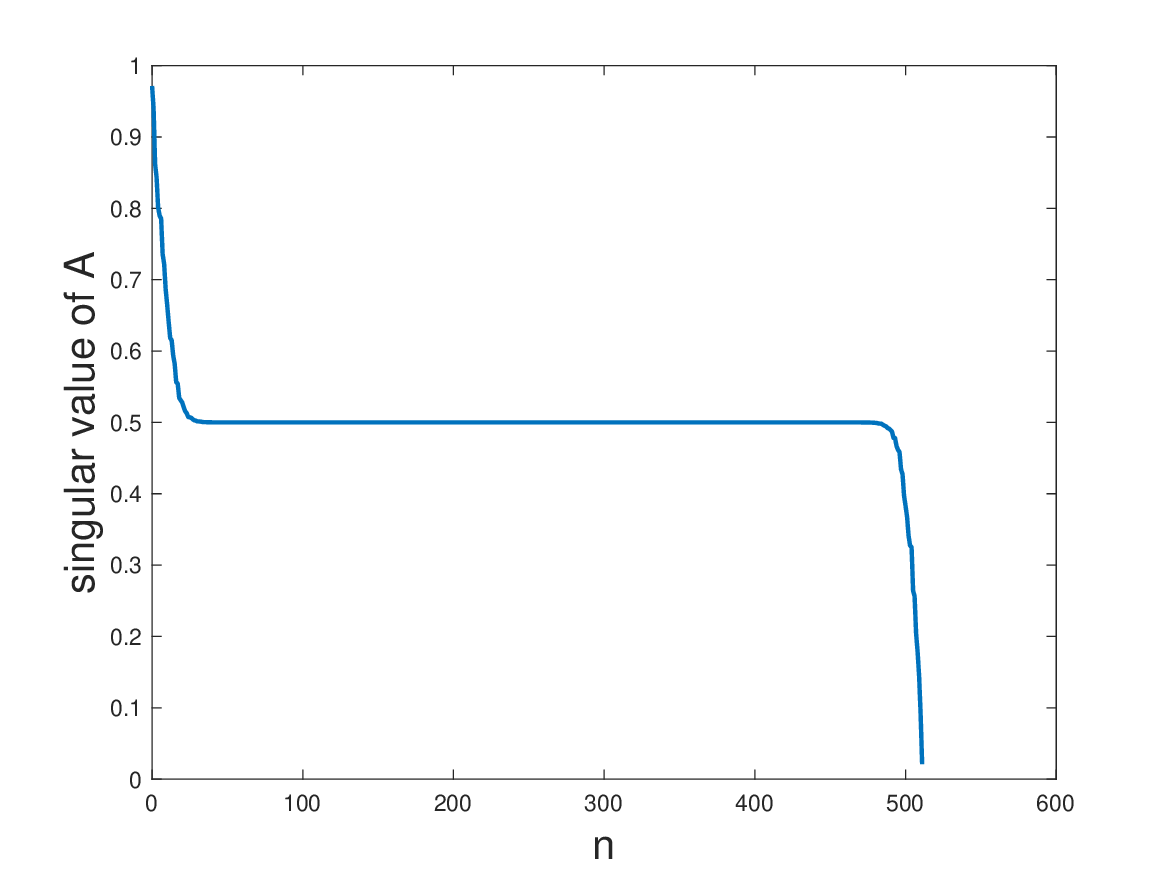, height=5cm} \epsfig{figure=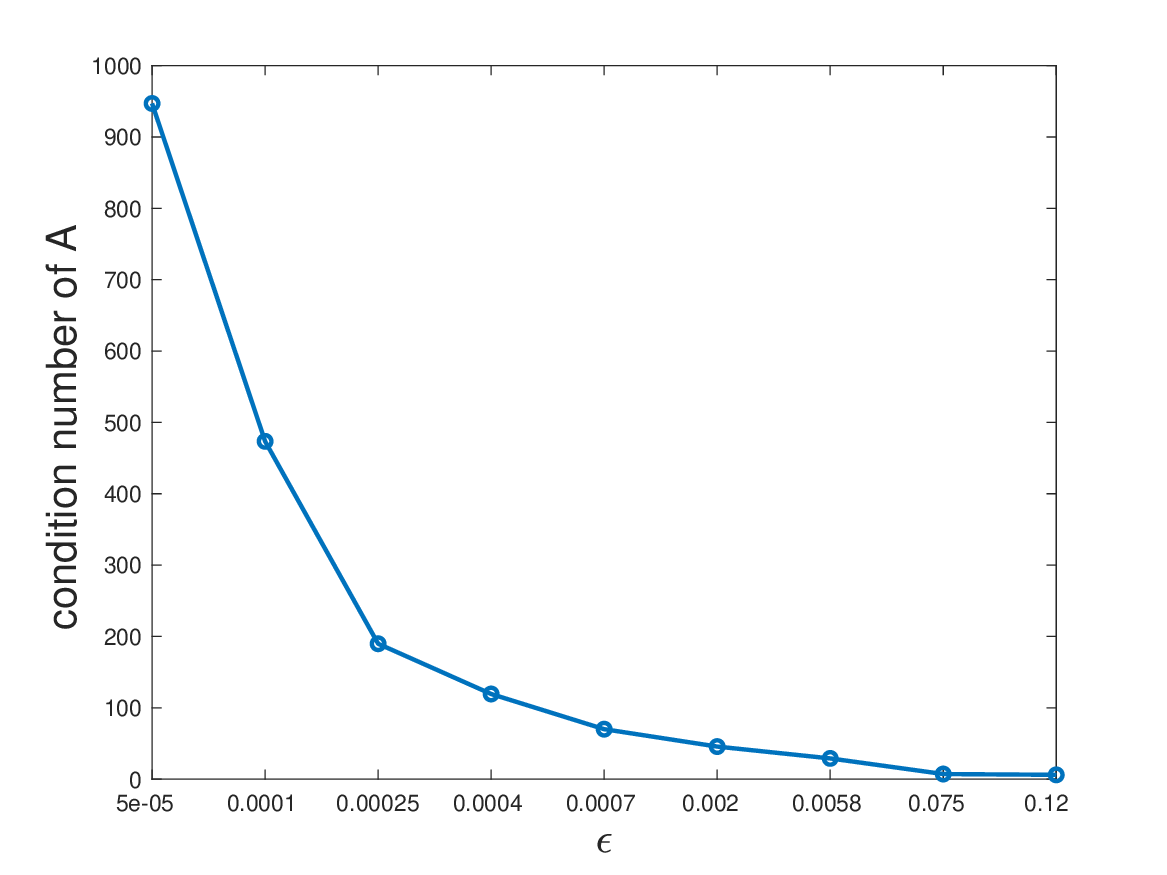, height=5cm}
\end{center}
\caption{Left: the singular values of $A$ in the decreasing order of $n$ when $\epsilon$ is $0.01$. Right: the condition numbers of $A$ as the distance $\epsilon$ tends to $0$. The dimension of $A$ is $512\times 512$.}
\label{singular}
\end{figure}

To solve $(\phi_1,\phi_2)$ from \eqref{sys-phi}, discretize each boundary $\p D_j$, $j=1,2$ into $N$ points, respectively. Let $x_j^k$, $k = 1,\dots,N$, be the nodal points on $\p D_j$. Then \eqref{sys-phi} becomes
$$
\left(-\frac{1}{2} I +A \right) \begin{bmatrix} \phi_1\\ \phi_2 \end{bmatrix} = \begin{bmatrix} Y_1\\ Y_2 \end{bmatrix},
$$
where 
\beq\label{A}
I = \begin{bmatrix} I_N & 0 \\ 0 & I_N \end{bmatrix}, \quad 
A = \begin{bmatrix} A_{11} & A_{12} \\
A_{21} &A_{22} \end{bmatrix}, 
\eeq
and
$$ Y_1 = - \begin{bmatrix} \frac{\p H}{\p \nu_{D_1}}(x_1^1) + \alpha_0 \frac{\p q_1}{\p \nu_{D_1}}\big|_{-} (x_1^1) \\ \vdots \\ \frac{\p H}{\p \nu_{D_1}}(x_1^N) + \alpha_0 \frac{\p q_1}{\p \nu_{D_1}}\big|_{-}(x_1^N) \end{bmatrix}, \quad  Y_2 = - \begin{bmatrix} \frac{\p H}{\p \nu_{D_2}}(x_2^1) + \alpha_0 \frac{\p q_2}{\p \nu_{D_2}}\big|_{-} (x_2^1)  \\ \vdots \\ \frac{\p H}{\p \nu_{D_2}}(x_2^N) + \alpha_0 \frac{\p q_2}{\p \nu_{D_2}}\big|_{-} (x_2^N) \end{bmatrix}.$$ 
Here $A$ is the evaluation of the kernel of $\mathbb{K}^*$. It is worth mentioning that the matrix $-\frac{1}{2} I +A$ has small singular values and the condition number of $A$ becomes worse as $\epsilon$ tends to zero, as shown in Figure \ref{singular}.
However, $\|\phi_1\|_{L^\infty(\p D_1)}$ and $\|\phi_2\|_{L^\infty(\p D_2)}$ are bounded regardless of $\epsilon$. 

\begin{figure}[!ht]
\begin{center}
\epsfig{figure=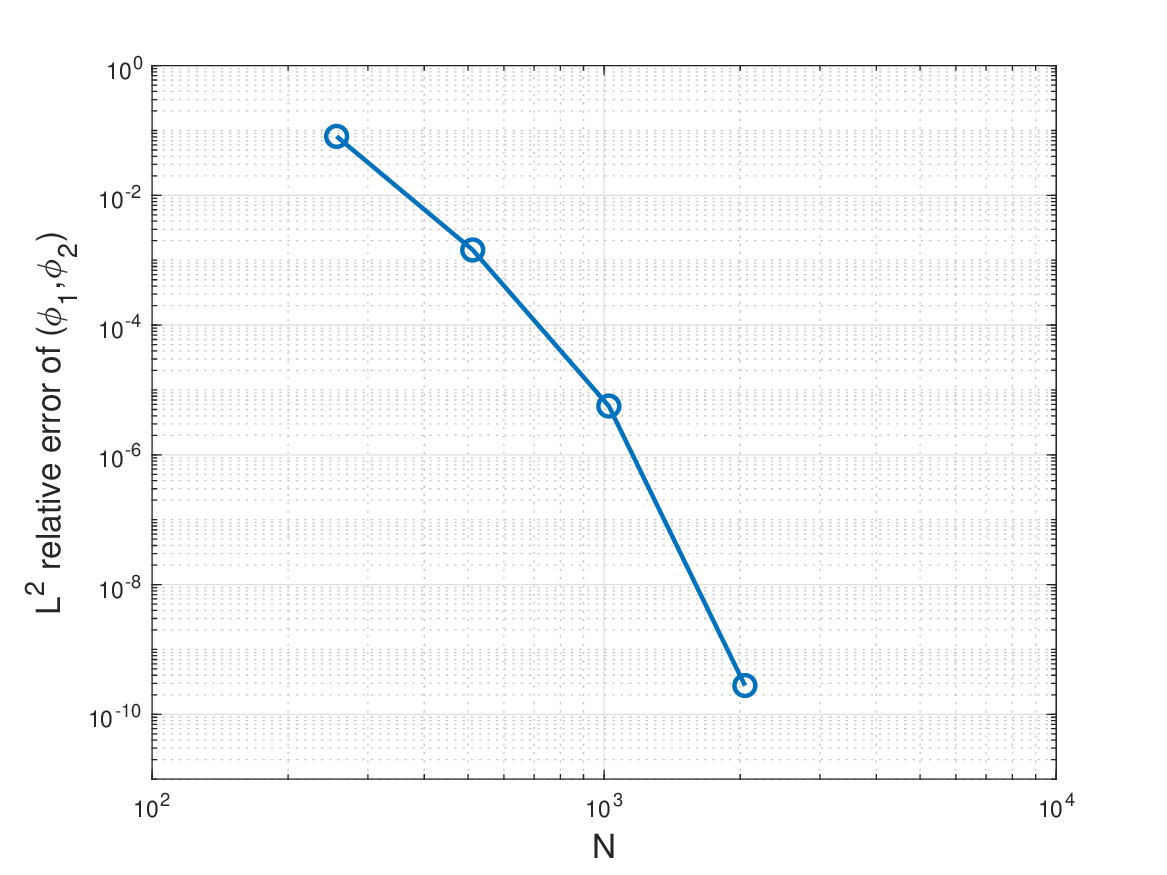, height=5cm} \epsfig{figure=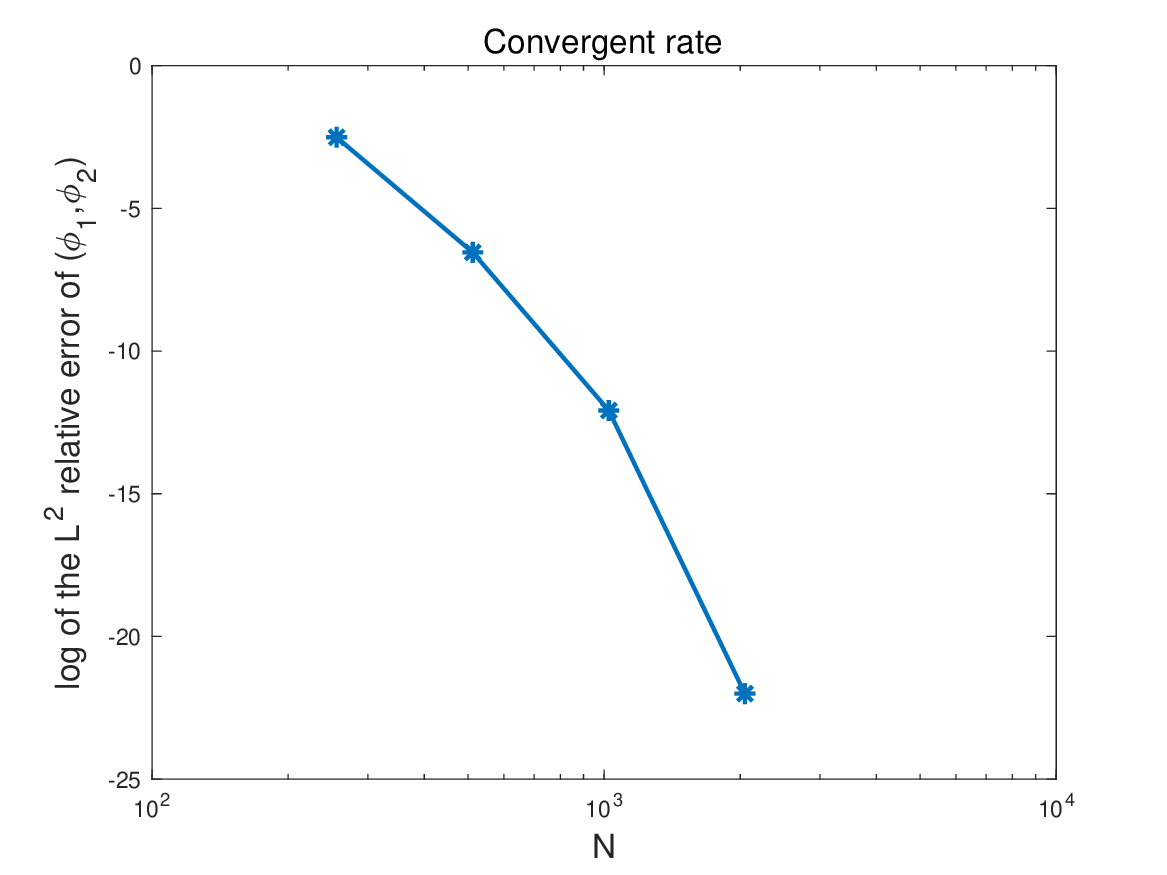, height=5cm}
\end{center}
\caption{Left: Relative error of $(\phi_1,\phi_2)$. Right: Logarithm of the relative error of $(\phi_1,\phi_2)$.}
\label{phi}
\end{figure}

We compute $(\phi_1,\phi_2)$ with $N = 256, 512, 1024, 2048, 4096$ equi-spaced points on $\p D_j$, $j = 1,2$, respectively. And then they are compared with the solution on the finer grid with $N = 4096$. Denote $(\phi_1^*,\phi_2^*)$ as the solution with grid number $N = 4096$. Let
$$\frac{ \|\phi_1 - \phi_1^* \|_{L^2(\p D_1)} }{ 2\|\phi_1^* \|_{L^2(\p D_1)}} + \frac{ \|\phi_2 - \phi_2^* \|_{L^2(\p D_2)} }{ 2\|\phi_2^* \|_{L^2(\p D_2)}}$$
be the relative $L^2$-errors of $(\phi_1,\phi_2)$ compared with $(\phi_1^*,\phi_2^*)$. Figure \ref{phi} (Left) shows that the relative errors decrease as the grid number $N$ increases. Figure \ref{phi} (Right) shows the logarithm of the relative error. One can see that $(\phi_1,\phi_2)$ converges very fast.

%%%%%%%%%%%%%%%%%%%%%%%%%%%%%%%%%%%%%%%%%%%
\subsection{Numerical scheme and effectiveness of the method}
%%%%%%%%%%%%%%%%%%%%%%%%%%%%%%%%%%%%%%%%%%

In this subsection we give the numerical scheme on the computation of the solution using characterization of the singular term method in the following Algorithm 1. We show the effectiveness of this method by comparing the convergent rate with the solution computed using layer potential techniques in a direct way.
\begin{algorithm}
\caption{Numerical scheme}
\begin{algorithmic}

\STATE \textbf{Step 1.}  Look for the solution to \eqref{main} in terms of the following form:
$$u(x) = \alpha_0 q(x) + H(x) + \mathcal{S}_{\p D_1}[\phi_1](x) + \mathcal{S}_{\p D_2}[\phi_2](x), \quad x\in \Rbb^2 \backslash \overline{(D_{1}\cup D_2)},$$
where $q$ is given by \eqref{q} and $\alpha_0$, $(\phi_1,\phi_2)$ are to be computed.

\STATE \textbf{Step 2.} Compute the stress concentrator factor $\alpha_0$ by considering the touching case \eqref{u-rho}:
\begin{itemize}
\item Discretize $\p D_\rho$ into $2N+N/8$ points;
\item Solve integral equation \eqref{psi} numerically;
\item Obtain $\alpha_0$ through \eqref{factor}.
\end{itemize}

\STATE \textbf{Step 3.} Compute the density functions $(\phi_1,\phi_2)$ through \eqref{sys-phi}:
\begin{itemize}
\item Discretize each $\p D_j$, $j = 1,2$ into $N$ points, respectively;
\item Compute the inward normal derivative $\frac{\p q_j}{\p \nu}\big|_{-}$, $j = 1,2$ by solving the Dirichlet problem \eqref{qj}, respectively;
\item Obtain $(\phi_1,\phi_2)$ by numerically solving \eqref{sys-phi}.
\end{itemize}

\STATE \textbf{Step 4.} Plot $u$.

\end{algorithmic}
\end{algorithm}

Denote $u^{deco}$ the solution computed following Algorithm 1. Denote $u^{dire}$ the solution to \eqref{main} by direct compution method. In fact, $u^{dire}$ can be written as
$$
u(x) = H(x)  + \mathcal{S}_{\p D_1}[\psi_1](x) + \mathcal{S}_{\p D_2}[\psi_2](x), \quad x\in \Rbb^2 \backslash {(D_{1}\cup D_2)},
$$
where two density functions $(\psi_1,\psi_2) \in L^2_0(\p D_1) \times L^2_0(\p D_2)$ satisfy 
$$
\left(-\frac{1}{2} I +A \right) \begin{bmatrix} \psi_1\\ \psi_2 \end{bmatrix} = -\begin{bmatrix} \frac{\p H}{\p \nu_{D_1}}\\[1ex] \frac{\p H}{\p \nu_{D_2}} \end{bmatrix},
$$
where $A$ is given by \eqref{A}. Note that the density functions $\psi_j$, $j = 1,2$ are as big as $1/\sqrt{\epsilon}$ near the origin point $(0,0)$ when $\nu\cdot \nabla H \neq 0$. Thus, applying single layer on $\psi_j$, $j = 1,2$, the error in the discretization of the single layer potential should become significant.

\begin{figure}[!ht]
\centering
\epsfig{figure=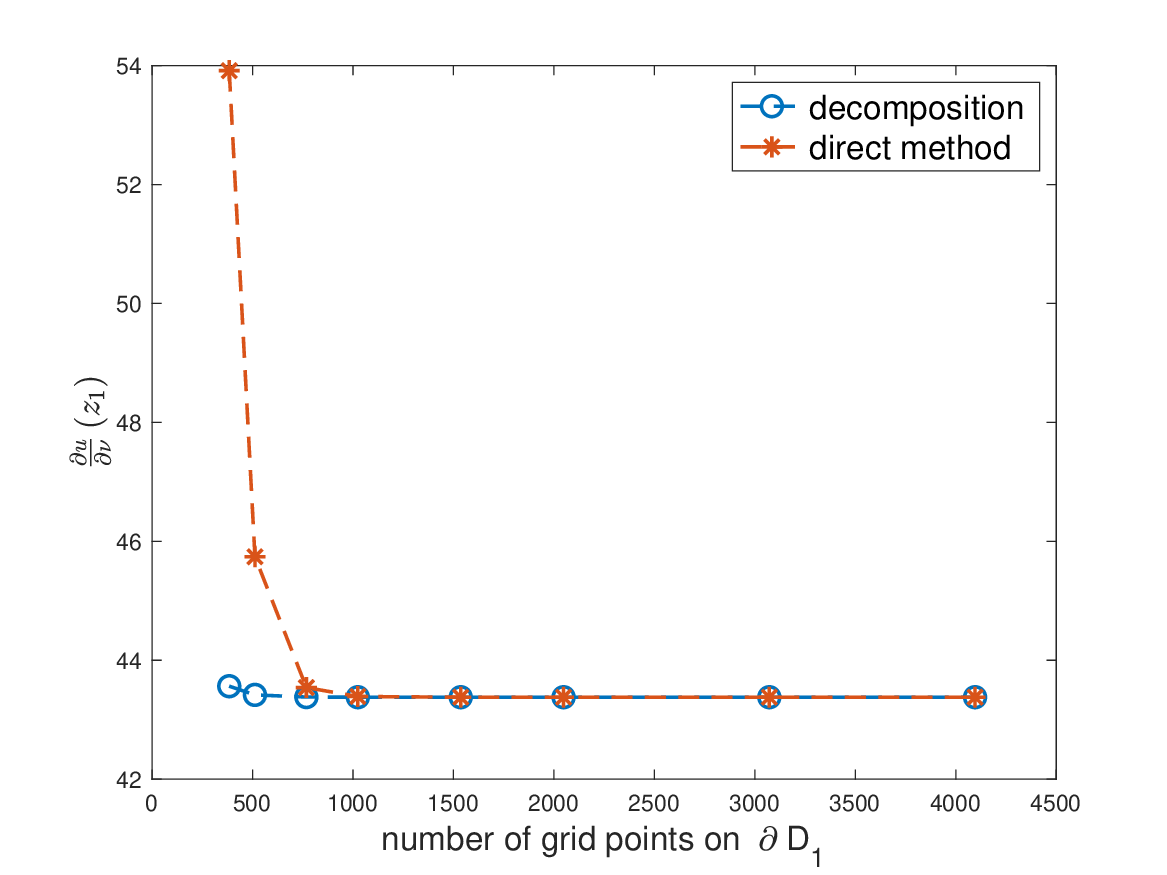, height=5cm} \epsfig{figure=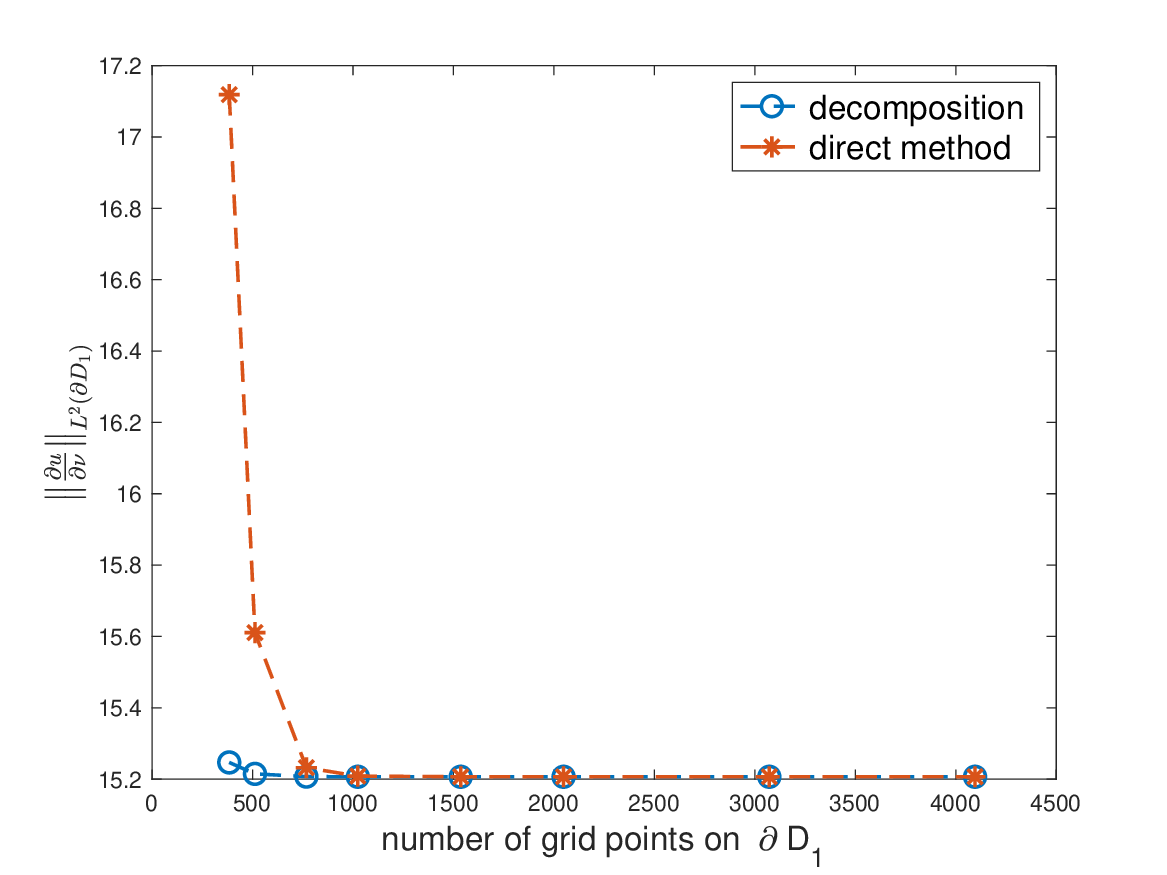, height=5cm}
\caption{Left: $\frac{\p u^{deco}}{\p \nu}$ and $\frac{\p u^{dire}}{\p \nu}$ at the closest point $z_1$; Right: $L^2$-norm of $\frac{\p u^{deco}}{\p \nu}$ and $\frac{\p u^{dire}}{\p \nu}$.}
\label{flux1}
\end{figure}
Let $D_1$ and $D_2$ be the same as in Section 3. Let the distance between $D_1$ and $D_2$ be $\epsilon = 0.01$. Let $z_1 = (-\epsilon/2,0)$ and $z_2 = (\epsilon/2,0)$ be the closest points on $\p D_1$ and $\p D_2$. The background potential is given by $H(x) = x_1$. Discretize each boundary $\p D_j$, $j = 1,2$ into $N$ points. Figure \ref{flux1} (Left) shows the normal derivative of $u^{dire}$(orange) and $u^{deco}$(blue) at $z_1$ for different values of $N$. Figure \ref{flux1} (Right) shows the $L^2$-norm of those on boundary $\p D_1$. One can see that both methods obtain convergent result, but the direct computation method needs finer meshes to obtain the accurate result. In fact, from Figure \ref{flux1} one can see that the direct computation method needs at least $N = 1024$ nodes on each boundary to obtain good result, whereas the characterization of the singular term method needs only needs $N = 512$. This further indicates that the scale of the line element of the discretization should be comparable with the scale of $\epsilon$ using the direct computation method, whereas the characterization of the singular term method does not limit to this restriction.

\begin{figure}[!ht]
\centering
\epsfig{figure=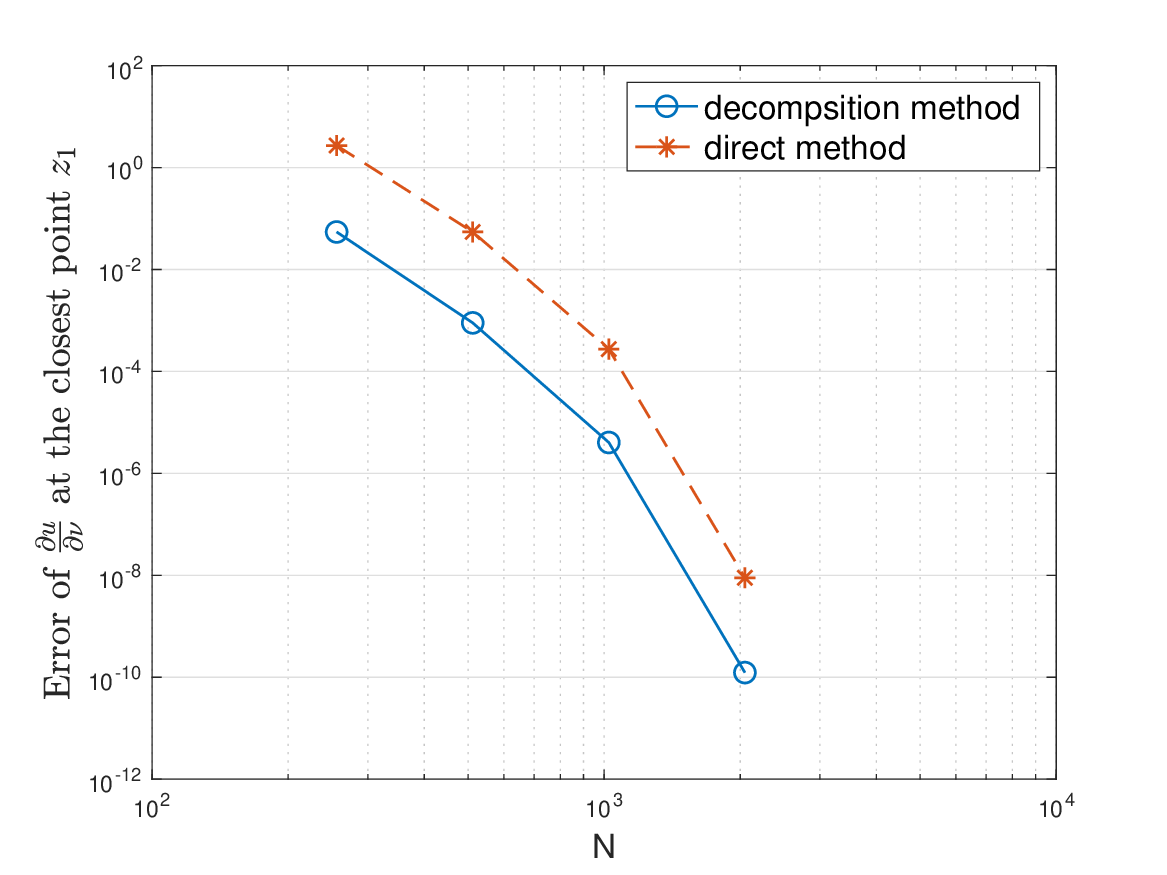, height=5cm} \epsfig{figure=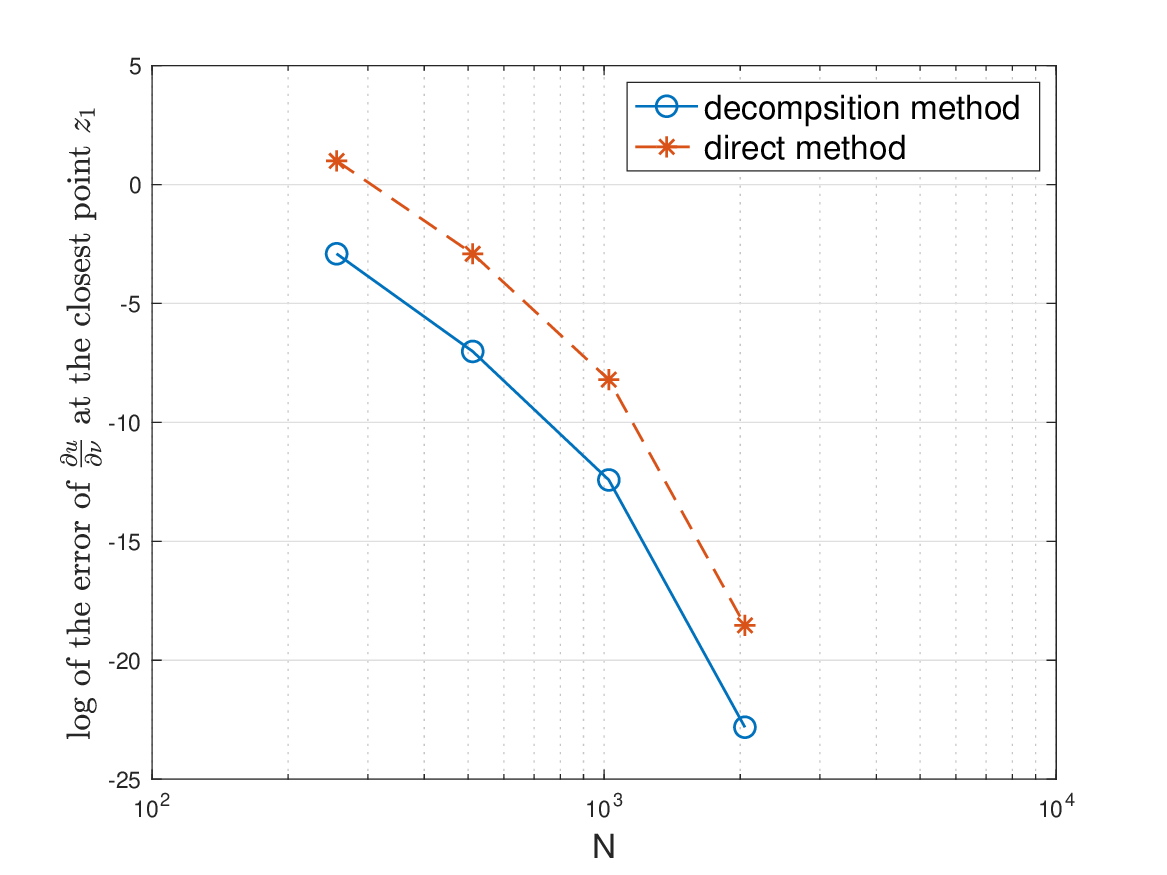, height=5cm}
\caption{Left: the relative errors of $\frac{\p u^{dire}}{\p\nu}$ and  $\frac{\p u^{deco}}{\p\nu}$ compared to $\frac{\partial u^{*}}{\partial \nu}$; Right: the logarithm of the relative errors. The distance $\epsilon=0.01$. The background potential is given by $H(x)=x_1$.}
\label{flux}
\end{figure}

From Figure \ref{flux1}, one can see that both methods obtain accurate results when the number of grid $N$ is sufficiently large. Hence, we assume that upon $N=4096$ the solution is regarded as the exact solution and denote it by $u^{*}$.

Fix $\epsilon = 0.01$. To show the effectiveness of the characterization of the singular term method, we compare the normal derivative $\frac{\p u^{dire}}{\p\nu}$ and  $\frac{\p u^{deco}}{\p\nu}$, with $\frac{\p u^{*}}{\p\nu}$, at the closest point $z_1$, for different values of $N = 256,512,1024,2048$. Let
$$\frac{\left| \frac{\p u^{dire}}{\p \nu_{D_1}}(z_1) - \frac{\p u^{*}}{\p \nu_{D_1}}(z_1) \right| } { \left|  \frac{\p u^{*}}{\p \nu_{D_1}}(z_1) \right| }  \quad \mbox{and} \quad  \frac{\left| \frac{\p u^{deco}}{\p \nu_{D_1}}(z_1) - \frac{\p u^{*}}{\p \nu_{D_1}}(z_1)\right| }{\left| \frac{\p u^{*}}{\p \nu_{D_1}}(z_1)\right|}  $$
be the relative errors of $\frac{\p u^{dire}}{\p\nu}$ and $\frac{\p u^{deco}}{\p\nu}$ at the closest point $z_1$, respectively. Figure \ref{flux} (Left) shows that the relative errors of $\frac{\p u^{dire}}{\p\nu}$ and $\frac{\p u^{deco}}{\p\nu}$ both decrease as the grid number $N$ increases. However, the error of $\frac{\p u^{deco}}{\p\nu}$ is much smaller than that of $\frac{\p u^{dire}}{\p\nu}$. The convergent speed of $\frac{\p u^{deco}}{\p\nu}$ is faster than that of $\frac{\p u^{dire}}{\p\nu}$, which is shown in Figure \ref{flux} (Right).

\begin{figure}[!ht]
\centering
\epsfig{figure=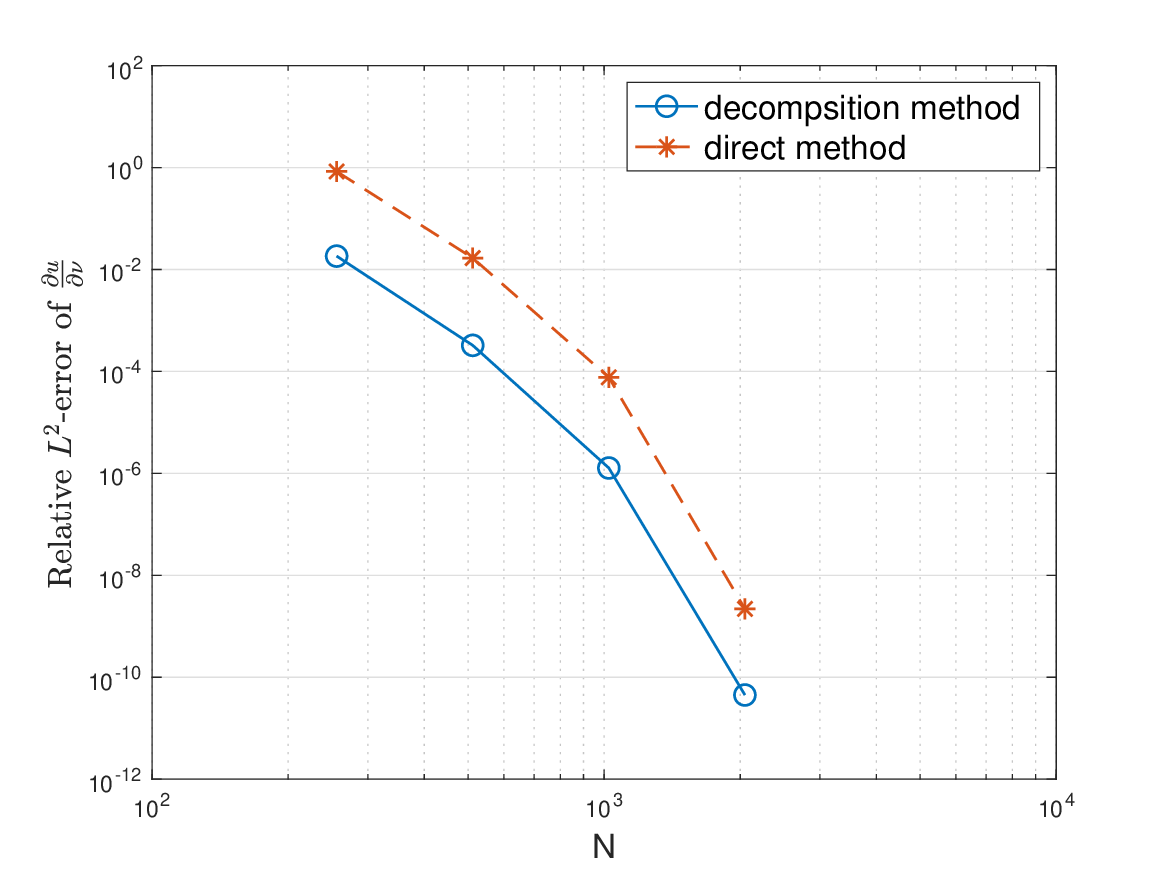, height=5cm} \epsfig{figure=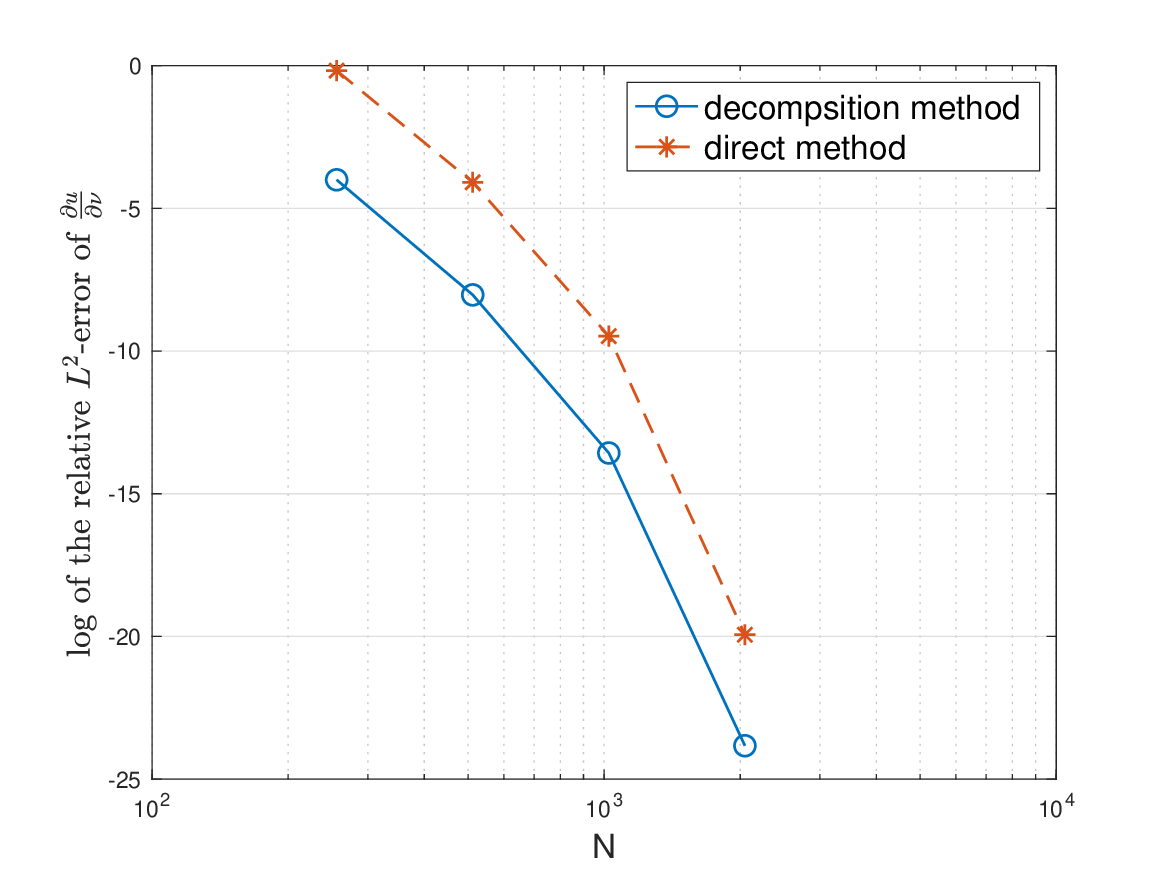, height=5cm}
\caption{Left: The relative $L^2$-errors of $\frac{\partial u^{dire}}{\partial \nu}$ and $\frac{\partial u^{deco}}{\partial \nu}$, respectively. Right: The logarithm of the relative $L^2$-error.}
\label{L2-error}
\end{figure}

We also compare the relative $L^2$-errors of $\frac{\p u^{dire}}{\p\nu}$ and $\frac{\p u^{deco}}{\p\nu}$, respectively, for different grid numbers $N = 256,512,1024,2048$. Let
$$ \frac{\|\frac{\p u^{dire}}{\p \nu_{D_1}} - \frac{\p u^{*}}{\p \nu_{D_1}}\|_{L^2(\p D_1)}}{2\|\frac{\p u^{*}}{\p \nu_{D_1}}\|_{L^2(\p D_1)}} + \frac{\|\frac{\p u^{dire}}{\p \nu_{D_2}} - \frac{\p u^{*}}{\p \nu_{D_2}}\|_{L^2(\p D_2)}}{2\|\frac{\p u^{*}}{\p \nu_{D_2}}\|_{L^2(\p D_2)}}$$
and
$$ \frac{\|\frac{\p u^{deco}}{\p \nu_{D_1}} - \frac{\p u^{*}}{\p \nu_{D_1}}\|_{L^2(\p D_1)}}{2\|\frac{\p u^{*}}{\p \nu_{D_1}}\|_{L^2(\p D_1)}} + \frac{\|\frac{\p u^{deco}}{\p \nu_{D_2}} - \frac{\p u^{*}}{\p \nu_{D_2}}\|_{L^2(\p D_2)}}{2\|\frac{\p u^{*}}{\p \nu_{D_2}}\|_{L^2(\p D_2)}}$$
be the relative $L^2$-errors of $\frac{\p u^{dire}}{\p\nu}$ and $\frac{\p u^{deco}}{\p\nu}$, respectively. Figure \ref{L2-error} (Left) shows that the relative $L^2$-errors of both methods decrease as the grid number $N$ increases, while the error of $\frac{\p u^{deco}}{\p\nu}$ is much smaller than that of $\frac{\p u^{dire}}{\p\nu}$. Figure \ref{L2-error} (Right) shows that the convergent rate of $\frac{\p u^{deco}}{\p\nu}$ is faster than that of $\frac{\p u^{dire}}{\p\nu}$, which indicates that the characterization of the singular term method is more effective. 

\begin{figure}[!ht]
\centering
\epsfig{figure=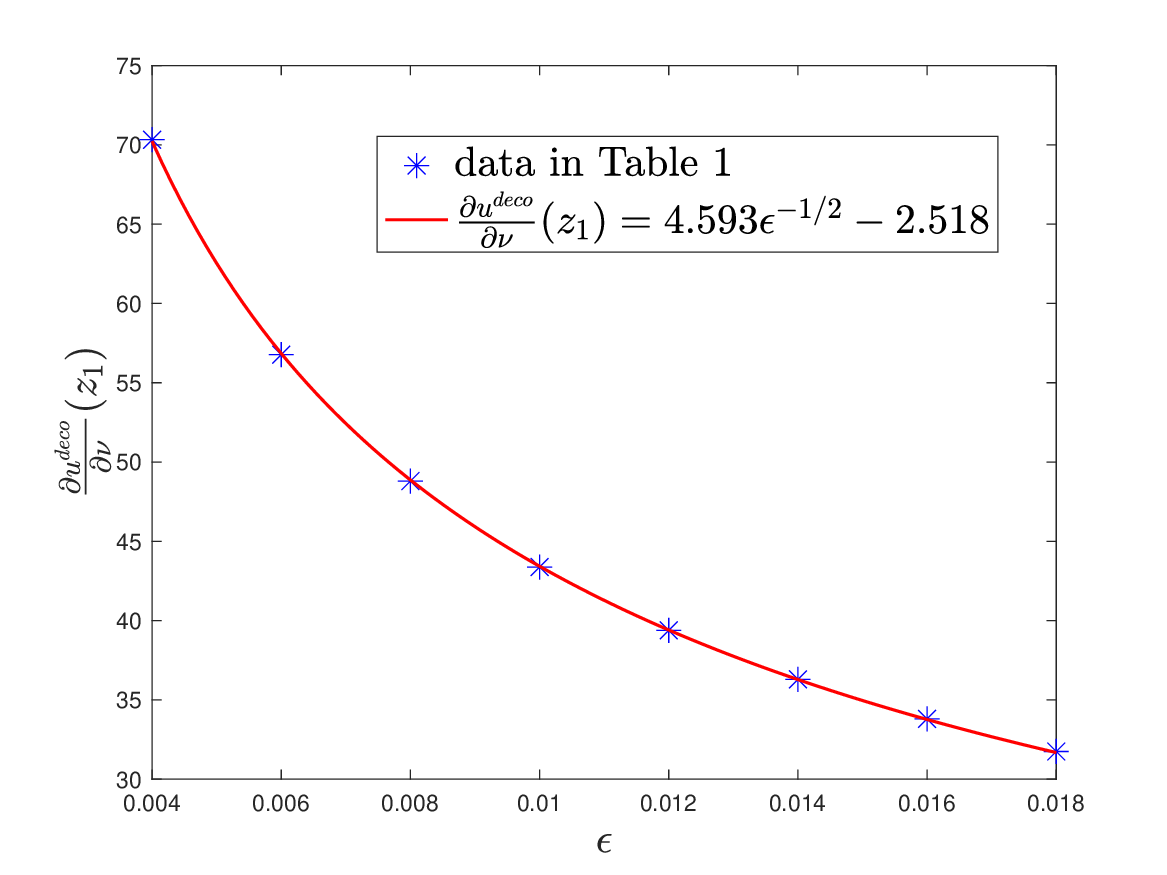, height=6cm}
\caption{Fitting curve and $\frac{\partial u^{deco}}{\partial \nu}$ at the closest point $z_1$ for different $\epsilon$.}
\label{eps}
\end{figure}

\begin{table}[!ht]
\centering
\begin{tabular}{ccc}
\hline
$\epsilon $  & $\frac{\partial u^{deco}}{\partial \nu}\left( z_1 \right)$ \\[1ex]\hline
0.018 & 31.745002  \\
0.016 & 33.802881  \\
0.014 & 36.292441  \\
0.012 & 39.387485  \\
0.010 & 43.378565  \\
0.008 & 48.798534  \\
0.006 & 56.770126  \\
0.004 & 70.326757  \\
\hline
\end{tabular}
\caption{The value of $\frac{\partial u^{deco}}{\partial \nu}$ at the closest point $z_1$, for different $\epsilon$.}
\label{data}
\end{table}
If we fix the number of grid points $N = 1024$ and vary the distance $\epsilon$ from $0.018$ to $0.004$. Table \ref{data} lists the values of the normal flux of $u^{deco}$ at the closest point $z_1$ for different values of $\epsilon$. The values are plotted as the blue star points in Figure \ref{eps}. The blow up rate of $\frac{\partial u^{deco}}{\partial \nu}$ is known to be $1/\sqrt{\epsilon}$ in two dimension. In fact, by \eqref{u-eps}, we have
$$\frac{\partial u^{deco}}{\partial \nu}(z_1) = \alpha_0 \frac{\partial q}{\partial \nu}(z_1) + O(1).$$
By \eqref{q} and \eqref{p12}, the above formula becomes
$$
\frac{\partial u^{deco}}{\partial \nu}(z_1) = \frac{\alpha_0}{2\pi} \left(\frac{z_1-p_1}{|z_1-p_1|^2} - \frac{z_1-p_2}{|z_1-p_2|^2} \right) + O(1).
$$
The concentration stress factor $\alpha_0$ is $10.312176$ which is given by Figure \ref{alpha_N} and Figure \ref{alpha_rho}. Together with the explicit value of $r_1 = 1/2$, $r_2 = 1/2$ in \eqref{p12}, we have
\beq\label{match}
\frac{\partial u^{deco}}{\partial \nu}(z_1) = 4.632\epsilon^{-1/2} + O(1).
\eeq
We now confirm the blow up rate $\epsilon^{-1/2}$ and its coefficient in \eqref{match} by fitting the values of $\frac{\partial u^{deco}}{\partial \nu}$ in Table \ref{data} with $\epsilon$ decreases from $0.018$ to $0.004$. The result is plotted as the red curve in Figure \ref{eps}. One can clearly see that the blow up rate is $\epsilon^{-1/2}$ and the coefficient of the fitting curve $4.593$ matches well with the coefficient in \eqref{match}. This result is interesting and reasonable. 

%%%%%%%%%%%%%%%%%%%%%%%%%%%%%%%%%%%%%%%%%%%%%%%%%%%%%%%%
\section{Numerical examples}
%%%%%%%%%%%%%%%%%%%%%%%%%%%%%%%%%%%%%%%%%%%%%%%%%%%%%%%%

In this section, we present some examples of numerical experiments on various different shapes of two closely located inclusions. The distance between two inclusions is $\epsilon = 0.01$.

\begin{figure}[!ht]
\centering
\epsfig{figure=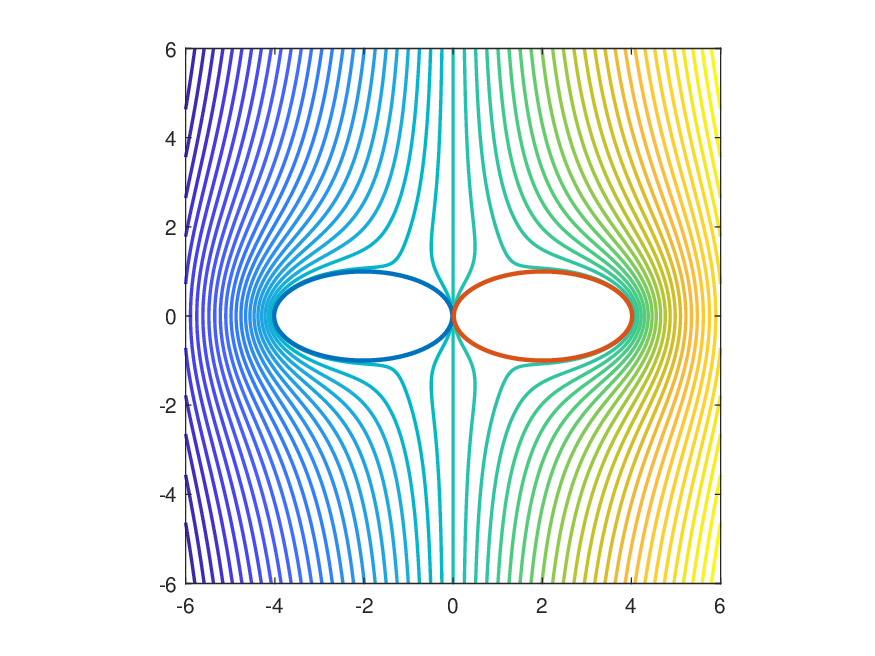, height=3.6cm}\epsfig{figure=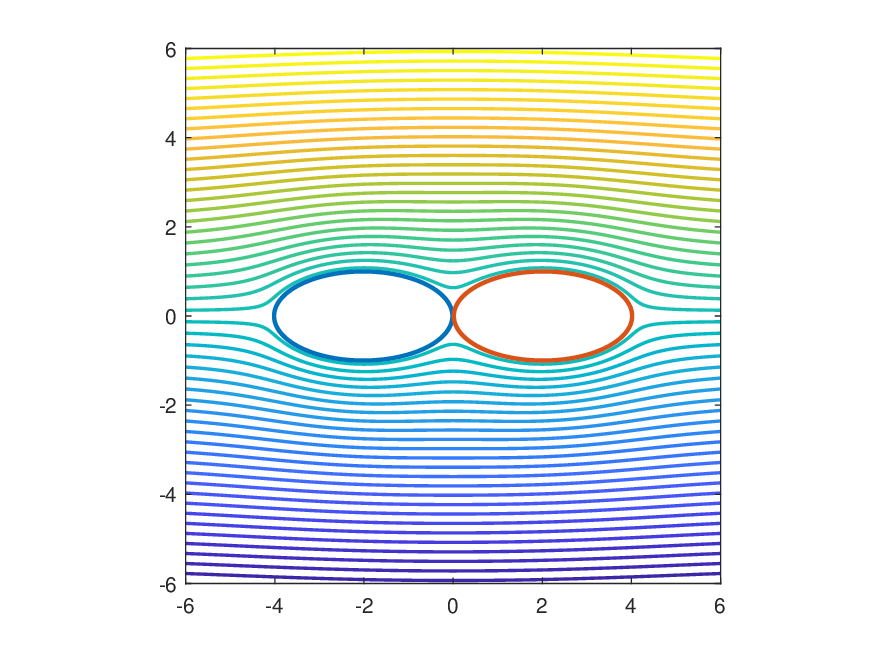, height=3.6cm}\epsfig{figure=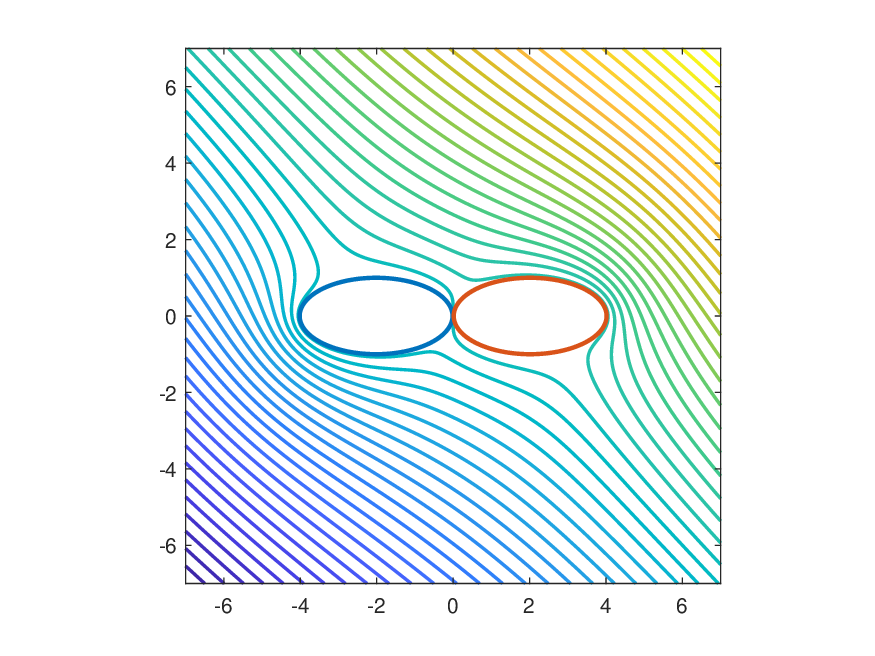, height=3.6cm}
\caption{Level curves of the stiff inclusions of elliptic shapes. Left: $H(x) = x_1$; Middle: $H(x) = x_2$; Right: $H(x) = x_1+x_2$.}
\label{ellipse}
\end{figure}

Firstly, let $D_1$ and $D_2$ be two elliptic inclusions of the same major axis $a=2$ and minor axis $b=1$, centered at $(-a-\epsilon/2,0)$ and $(a+\epsilon/2,0)$, respectively. Discretize each boundary $\p D_j$, $j =1,2$, into 256 grid nodes. Applying linear field $H(x) = x_1$, Figure \ref{ellipse} (Left) shows the uniformly spaced contour level curves. Applying linear field $H(x) = x_2$, one can see from Figure \ref{ellipse} (Middle) that the gradient does not blow up. Let $H(x) = x_1+x_2$, Figure \ref{ellipse} (Right) shows the uniformly spaced contour level curves.

\begin{figure}[!ht]
\centering
\epsfig{figure=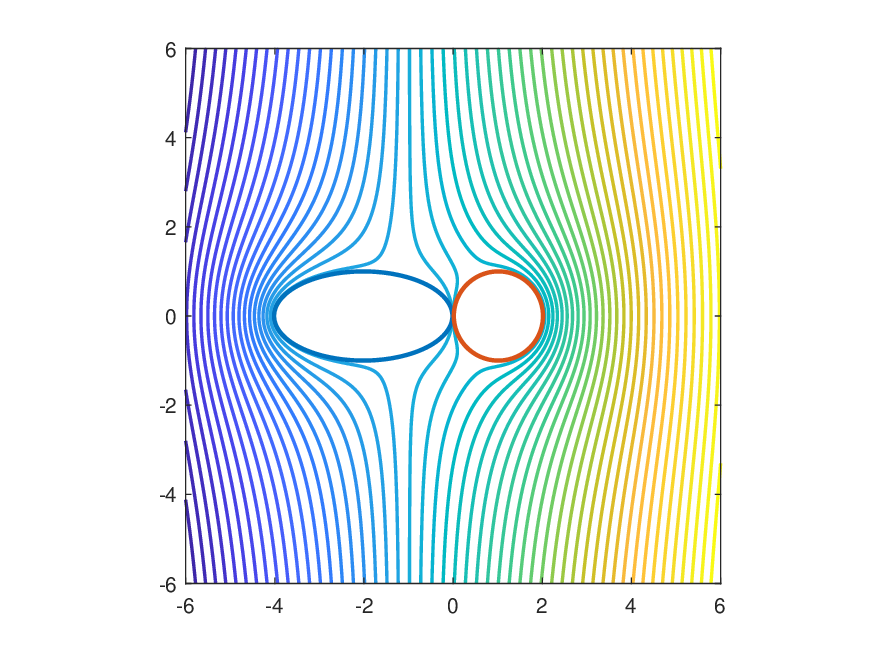, height=3.6cm}\epsfig{figure=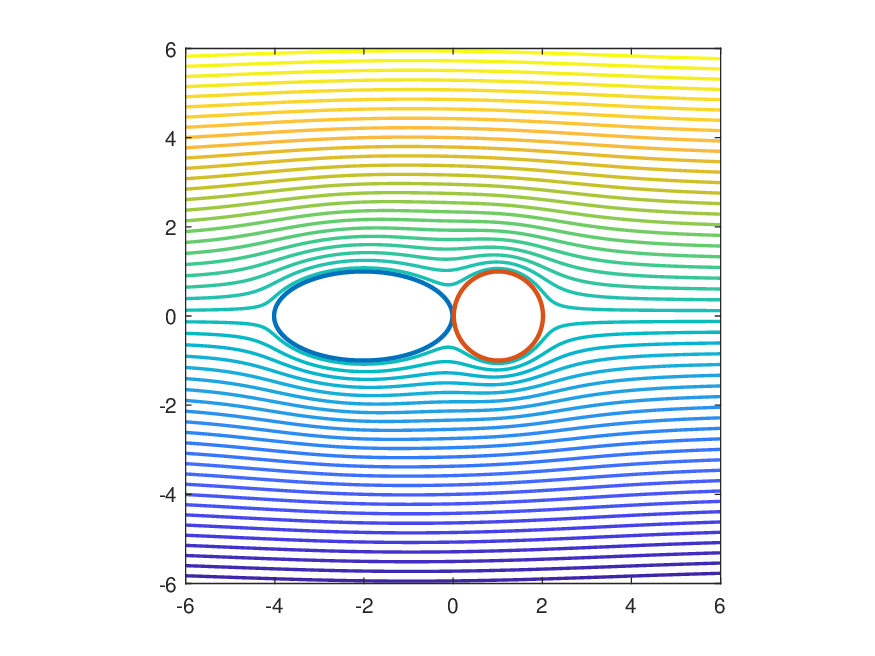, height=3.6cm}\epsfig{figure=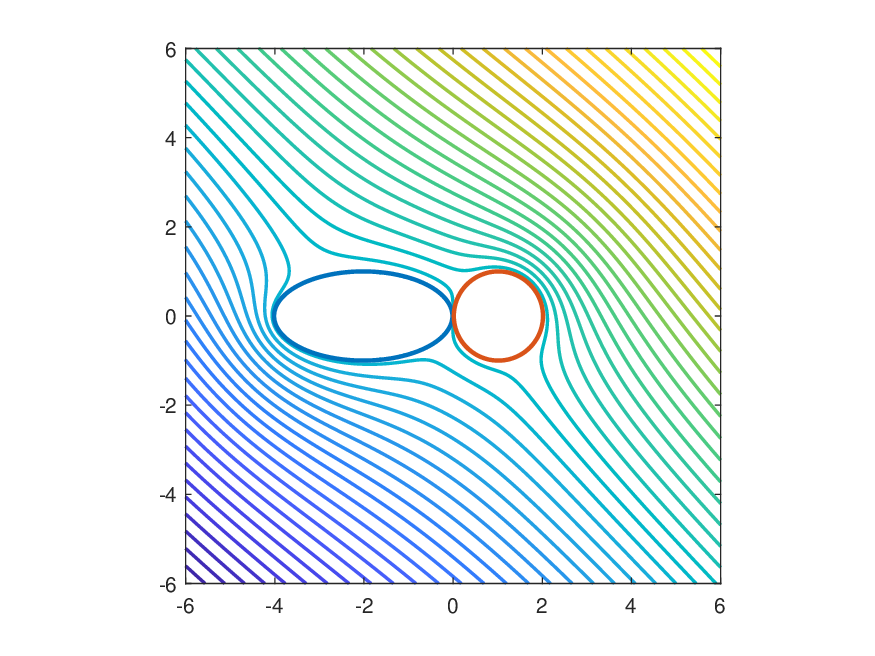, height=3.6cm}
\caption{Level curves of the stiff inclusions of elliptic and disk shapes. Left: $H(x) = x_1$; Middle: $H(x) = x_2$; Right: $H(x) = x_1+x_2$.}
\label{ellipse-disk}
\end{figure}

Secondly, let $D_1$ be an ellipse with the major axis $a = 2$ and minor axis $b = 1$, centered at $(-a-\epsilon/2,0)$, and let $D_2$ be a circle of radius $r=1$ centered at $(r+\epsilon/2,0)$. Figure \ref{ellipse-disk} shows the uniformly spaced contour level curves when $H(x) = x_1$, $H(x) = x_2$ and $H(x) = x_1+x_2$, respectively.

\begin{figure}[!ht]
\centering
\epsfig{figure=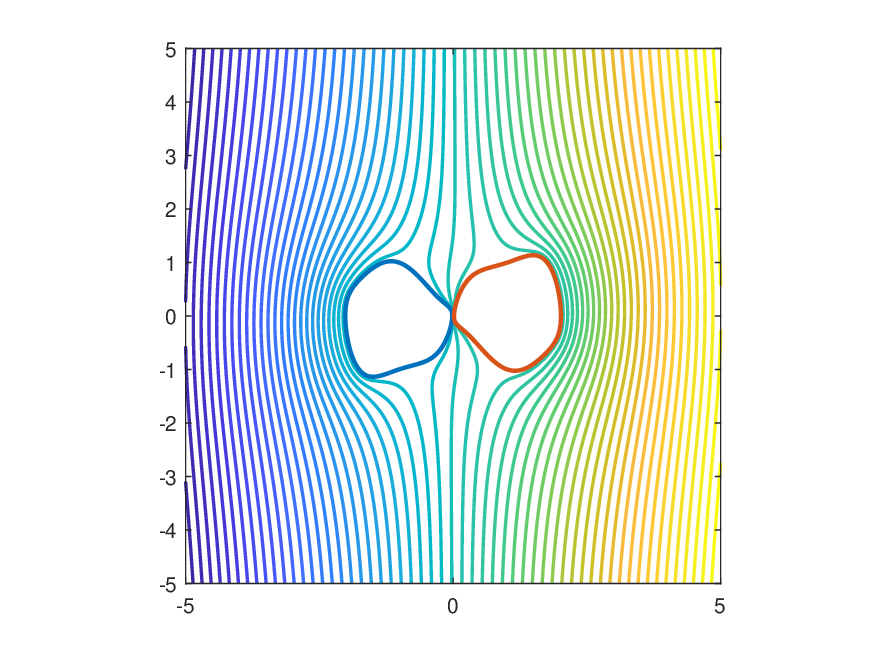, height=3.6cm}\epsfig{figure=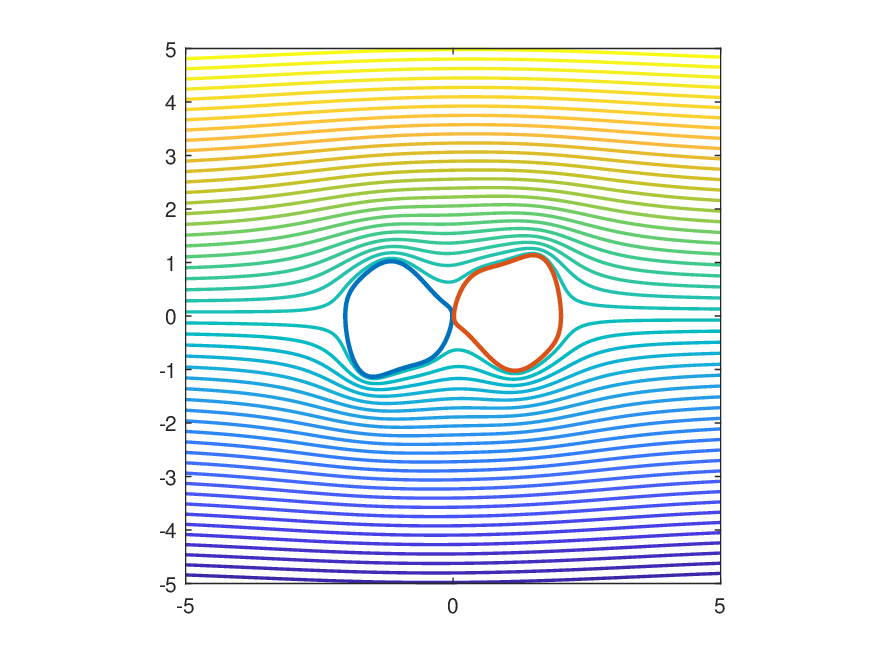, height=3.6cm}\epsfig{figure=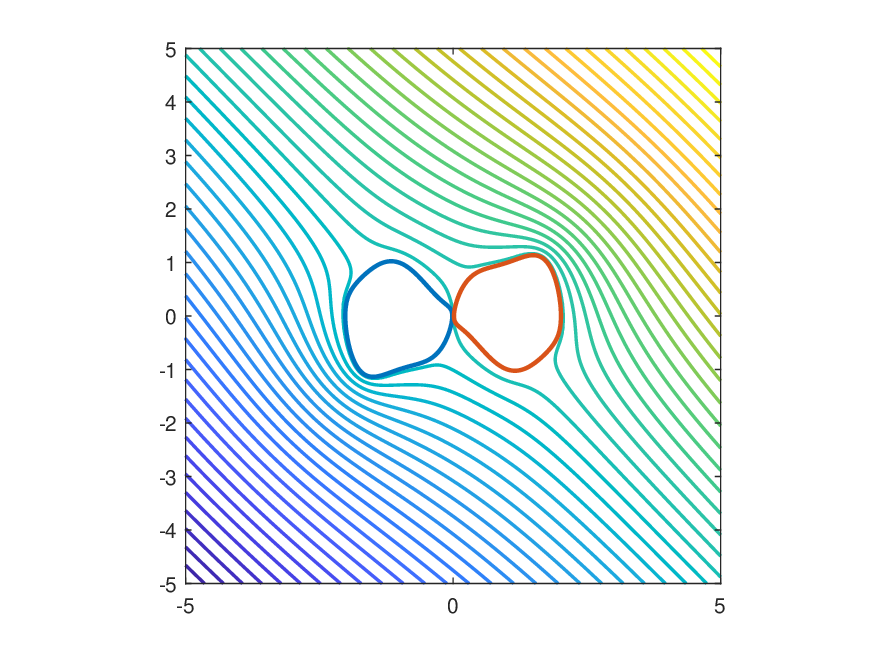, height=3.6cm}
\caption{Level curves of the stiff inclusions of general shape. Left: $H(x) = x_1$; Middle: $H(x) = x_2$; Right: $H(x) = x_1+x_2$.}
\label{general}
\end{figure}

As the final example, Figure \ref{general} shows the uniformly spaced contour level curves for two stiff inclusions of general shape, when $H(x) = x_1$, $H(x) = x_2$ and $H(x) = x_1+x_2$, respectively. The boundaries of two inclusions are given by the following parametrization functions for $\theta \in [0,2\pi)$:
$$
\begin{cases}
x_1 = -\frac{\epsilon}{2} - 1 + \cos(\theta), \\
x_2 = -\frac{1}{12} + \sin(\theta) - \frac{1}{6} \sin(2\theta) + \frac{1}{12} \cos(4\theta).
\end{cases}
$$

%%%%%%%%%%%%%%%%%%%%%%%%%%%%%%%%%%%%%%%%%%%%%%%%%%%%%%%%
\section{Conclusion}
%%%%%%%%%%%%%%%%%%%%%%%%%%%%%%%%%%%%%%%%%%%%%%%%%%%%%%%%

In this paper, we show through numerical simulations that the computation of the stress concentration between closely located stiff inclusions of general shapes can be realized by only using regular meshes. Using the characterization of the singular term method, we can decompose the solution into a singular and a regular term. After extracting the singular in a precise way, we can compute the remaining term using regular meshes. The key point in our computation lies in the computation of the stress concentration factor as well as the singular term. We have shown that the computation of the stress concentration factor converges very fast. By comparing the convergent rate with the solution computed using layer potential techniques in a direct way, we conclude that the characterization of the singular term method can be used effectively for the computation of the solution. This numerical method could also be generalized to the three dimensional case, which we will investigate in a forthcoming work.

%%%%%%%%%%%%%%%%%%

%\section*{Acknowledgement}
%The authors would like to express their gratitude to Hyeonbae Kang for his discussion on this  work. The work of the authors was supported by the NSF of China grant No. 11901523.

%%%%%%%%%%%%%%%%%%%%%%%%%%%%%%%%%%%%%%%%%%%%%%%

\end{document}